\documentclass[10pt, oneside, a4paper]{article}

\usepackage{amsmath,amssymb,amsfonts}
\usepackage{dsfont}
\usepackage{graphicx}
\usepackage{amsthm}
\usepackage{hyperref}

\newcommand{\N}{\mathbb{N}}
\newcommand{\Z}{\mathbb{Z}}
\newcommand{\R}{\mathbb{R}}
\newcommand{\tm}{\times}
\newcommand{\ep}{\varepsilon}
\newcommand{\unit}{\mathds{1}}
\newcommand{\id}{\mathrm{id}}
\newcommand{\inner}{\mathrm{int}}
\newcommand{\Fix}{\mathrm{Fix}}
\newcommand{\rme}{\mathrm{e}}

\newcommand{\FC}{\mathcal{F}}
\newcommand{\GC}{\mathcal{G}}
\newcommand{\IC}{\mathcal{I}}
\newcommand{\JC}{\mathcal{J}}
\newcommand{\KC}{\mathcal{K}}
\newcommand{\LC}{\mathcal{L}}
\newcommand{\NC}{\mathcal{N}}

\newcommand{\bsign}[1]{\bar{#1}}

\newtheorem{proposition}{Proposition}[section]
\newtheorem{definition}{Definition}[section]
\newtheorem{lemma}{Lemma}[section]
\newtheorem{theorem}{Theorem}[section]
\newtheorem{example}{Example}[section]
\newtheorem{remark}{Remark}[section]
\newtheorem{corollary}{Corollary}[section]

\def\keywords{\vspace{.5em}
{\bf{\textit{{\small Keywords}}--}\,\relax%
}}

\allowdisplaybreaks
\begin{document}

\title{The Small-Gain Condition for Infinite Networks Modeled on $\ell^{\infty}$-Spaces}
\author{Christoph Kawan\footnote{Ulm, Germany; e-mail: \href{mailto:christoph@kawan@gmx.de}{christoph.kawan@gmx.de}}\,\,\footnote{Acknowledgements: The author thanks Andrii Mironchenko for many useful comments on previous versions of this manuscript and pointers to the literature.}}

\date{}
\maketitle

\begin{abstract}
In recent years, attempts have been made to extend nonlinear small-gain theorems for input-to-state stability (ISS) from finite networks to countably infinite networks with finite indegrees. Under specific assumptions about the interconnection gains and the ISS formulation, corresponding infinite-dimensional small-gain results have been proven. However, concerning these assumptions, the results are still too narrow to be considered a full extension of the state-of-the-art for finite networks. We take a step to closing this gap by developing a general technical framework within which the small-gain condition for both finite and infinite networks can be analyzed. This includes a thorough investigation of various monotone operators associated with a network and a specific ISS formulation. Our results extend and generalize the existing theory for finite networks and yield complete characterizations of the small-gain condition for specific ISS formulations. They also provide much deeper insight into the different characterizations of the small-gain condition, including the no-joint-increase condition $\Gamma(s) \not\geq s$.
\end{abstract}


\section{Introduction}

Given a stabilizing feedback controller for a nonlinear system, it is possible that small deviations from this controller lead to divergent behavior. Since any practical implementation of a controller necessarily deviates from the conceptual controller, this is a serious problem. As a remedy, Eduardo Sontag introduced the concept of input-to-state stability (ISS) in \cite{Son}. If the open-loop system obtained by including time-dependent bounded perturbations into the feedback loop is input-to-state-stable, practical stability is guaranteed and the bounded deviation of trajectories from the equilibrium can be estimated quantitatively over an infinite time horizon. As any general property of a system can be hard to verify in a concrete example, it is necessary to find sufficient conditions for ISS which can be checked by computation. Sontag also introduced the notion of an ISS Lyapunov function, whose existence is equivalent to ISS, at least in finite dimensions. This shifts the problem of verifying ISS directly to the construction (or proof of existence) of an ISS Lyapunov function, which is not necessarily easier though, especially for high-dimensional systems.%

The small-gain approach tries to mitigate the difficulty by splitting the task of finding a Lyapunov function into two manageable subtasks. The idea is to decompose the system into a number of lower-dimensional, physically interconnected subsystems. Every subsystem has two types of inputs, the external perturbations acting on the overall system and the outputs of the other subsystems. If one can find an ISS Lyapunov function for each subsystem with respect to both types of inputs, there is a chance that these functions can be composed to an overall ISS Lyapunov function. This strategy works if the interconnections are sufficiently weak, which is formalized in terms of a so-called \emph{small-gain condition}. The latter can be formulated as a property of a monotone operator $\Gamma$ (called \emph{gain operator}), built from the interconnection gains and acting on the positive cone of an ordered vector space. In finite dimensions, there are four equivalent ways to describe the small-gain condition: in terms of a path evolving in the decay set of $\Gamma$, in terms of an asymptotic stability property of the dynamical system generated by $\Gamma$, in terms of a monotone bounded invertibility property of the operator $\id - \Gamma$, and in terms of the componentwise action of $\Gamma$ in one time step. The last characterization is particularly useful, since it allows for computational verification in many cases. There is also a connection to the theory of positive linear operators. Indeed, under specific assumptions about the gains and the ISS formulation, $\Gamma$ becomes a positive linear operator and the small-gain condition is equivalent to the spectral radius condition $r(\Gamma) < 1$.%

An alternative view on the small-gain theory is that the network structure is given a priori by the setup of the control system, which naturally splits into many interacting components that one wants to control simultaneously, and which may be distributed over a large area. Examples can be found, e.g., in automated traffic systems, electric power grids, or aerospace and defense systems. In such applications, the number of subsystems can vary in time and become quite large. As the number of subsystems in a network grows, measures of stability may deteriorate even if ISS is preserved. A successful approach to dealing with large systems often is to over-approximate them by some sort of infinite systems. This strategy is used, for instance, in physics to model the dynamics of a fluid or a gas, consisting of an enormous number of interacting molecules. Motivated by such ideas, we study the ISS small-gain approach for networks consisting of an infinite number of interacting systems. First results in this direction have been obtained in \cite{Dea,DPa,Kea,KMZ,KZa,KMG,Mea}. However, in all of these papers, the authors either restricted their attention to a narrow class of gains and/or formulations of the ISS property, or did not obtain a full generalization of the existing results for finite networks as presented in \cite{DR1,DR2,KJi,Ruf}.%

In this paper, we aim at such a generalization. To this end, we introduce a very general class of gain operators associated with a network of countably many interconnected systems. We mainly impose one restriction, namely that each subsystem is only influenced by finitely many other subsystems. For such a setup, the paper \cite{KMZ} already provided a small-gain theorem in which the small-gain condition was formulated in terms of the existence of a so-called \emph{path of strict decay}. However, the existence of such a path was only investigated for a particular formulation of the ISS property, leading to gain operators with a very convenient algebraic property. In this paper, we aspire much greater generality.%

We abstain from an extensive literature review on the control of infinite networks or input-to-state stability of infinite-dimensional systems, since one can find such reviews in the papers cited above. We just mention \cite{Bea,Be2} as two prominent examples, which have connections to the topic of this paper, but address control problems for large-scale or infinite networks rather than mere stability problems. Our paper also abstains from the discussion of concrete examples, since they would not provide a valuable contribution to the general theory without a much deeper analysis tailored to those examples, which can be part of future work. Finally, we note that our results can potentially contribute to the general theory of monotone operators on Banach spaces, but this needs some further abstraction, which can also be part of future work.%

We now describe the structure of the paper and the contents of each section.%

In Section \ref{sec_gain_ops}, we define gain operators together with the main object needed to construct an overall ISS Lyapunov function from ISS Lyapunov functions of the subsystems, namely a \emph{path of strict decay}.%

{\bf Subsection \ref{subsec_prelim}}: We introduce the fundamental mathematical objects and concepts with the corresponding notation, in particular the ordered Banach space $(\ell^{\infty}(\IC),\ell^{\infty}_+(\IC))$ for a countable (finite or infinite) index set $\IC$ with the induced partial order $\leq$ on $\ell^{\infty}(\IC)$, the binary operation $\oplus$ on $\ell^{\infty}_+(\IC)$ (componentwise maximum), the vector $\unit = (1,1,1,\ldots)$, and discrete-time dynamical systems on $\ell^{\infty}_+(\IC)$ generated by monotone operators.%

{\bf Subsection \ref{subsec_def_props}}: We define gain operators and introduce the notion of a path of strict decay. A gain operator $\Gamma$ acts as a monotone operator on the cone $\ell^{\infty}_+(\IC)$. A path $\sigma$ of strict decay evolves in the decay set $\Psi(\Gamma_{\rho}) := \{ s \in \ell^{\infty}_+(\IC) : \Gamma_{\rho}(s) \leq s \}$ of an \emph{enlarged gain operator} $\Gamma_{\rho} = (\id + \rho) \circ \Gamma$, where $\rho$ is a $\KC_{\infty}$-function. By definition, such a path needs to satisfy a bunch of properties, including a local bi-Lipschitz condition. We prove that some of these properties are redundant. Under considerably weaker assumptions on a path $\sigma$, it is possible to obtain the desired path by certain algebraic manipulations performed on $\sigma$. This leads to the notion of a \emph{$C^0$-path (of strict decay)}, which is a continuous and increasing path in $\Psi(\Gamma_{\rho})$ satisfying coercivity estimates.%

{\bf Subsection \ref{subsec_associated_ops}}: We introduce further monotone operators derived from a gain operator together with the associated dynamical systems. This allows us to characterize the existence of a path with all the properties of a $C^0$-path except for continuity.%

{\bf Subsection \ref{subsec_stability}}: We investigate stability properties of the dynamical system $\Sigma(\Gamma)$ generated by a gain operator $\Gamma$. We first characterize global attractivity of $\Sigma(\Gamma)$ in the weak$^*$-topology. Here, the classical small-gain condition $\Gamma(s) \not\geq s$ for $s > 0$, also called \emph{no-joint-increase (NJI) condition}, is an essential ingredient. We then prove that global attractivity in the norm implies a uniform NJI condition, requiring that for each index $i \in \IC$, the componentwise decay condition $\Gamma_j(s) < s_j$ holds for some index $j = j(i)$, uniformly close to $i$ in the topology of the interconnection graph. We call this condition the \emph{topologically uniform NJI condition}. For finite networks, this condition is equivalent to the NJI condition, and in this sense, it is a proper generalization to the infinite case. We also prove a converse result, showing that the topologically uniform NJI condition for $\Gamma_{\rho}$ together with the cofinality of $\Psi(\Gamma_{\rho})$ in $\ell^{\infty}_+(\IC)$ implies global attractivity.%

In Section \ref{sec_path_of_decay}, we investigate conditions for the existence of a $C^0$-path.%

{\bf Subsection \ref{subsec_max_mbi}}: We introduce a property of gain operators, called the \emph{maximum monotone bounded invertibility ($\oplus$-MBI) property}. This property, in particular, implies the uniform cofinality of $\Psi(\Gamma)$ and weak$^*$ global attractivity of $\Sigma(\Gamma)$. Moreover, it implies that each of the so-called \emph{projected gain operators} $\bsign{\Gamma}_b(s) := b \oplus \Gamma(s)$ with $b \in \ell^{\infty}_+(\IC)$ has a minimal and a maximal fixed point $s_*(b)$ and $s^*(b)$, respectively, such that the order interval $[s_*(b),s^*(b)]$ is a compact global attractor of $\Sigma(\bsign{\Gamma}_b)$ in the weak$^*$-topology. Finally, the $\oplus$-MBI property admits a characterization in terms of another uniform NJI condition (called the \emph{spatially uniform NJI condition}).%

{\bf Subsection \ref{subsec_necessary}}: We obtain as necessary conditions for the existence of a $C^0$-path (i) the $\oplus$-MBI property for an enlarged gain operator $\Gamma_{\rho}$ and (ii) \emph{uniform global asymptotic stability (UGAS)} of $\Sigma(\Gamma_{\rho})$.%

{\bf Subsection \ref{subsec_path_construction}}: We introduce a method to construct a $C^0$-path, which uses the maximal fixed points $s^*(r\unit)$ in combination with local linear interpolation between the points of $\Sigma(\bsign{\Gamma}_{r\unit})$-trajectories. Using this method, we establish that the $\oplus$-MBI property is equivalent to the existence of a path, satisfying all properties of a $C^0$-path except for norm-continuity from above, which has to be replaced by sequential weak$^*$-continuity. When the index set $\IC$ is finite, a complete characterization of the existence of a path of strict decay in terms of the $\oplus$-MBI property follows. Moreover, we conclude that the $\oplus$-MBI property for $\Gamma_{\rho}$ implies that each finite subnetwork admits a path of strict decay with coercivity bounds that are independent of the chosen subnetwork.%

In Section \ref{sec_final_theorem}, we formulate our main results, which combine the partial results described above. In particular, we see that the existence of a path of strict decay implies both uniform $\oplus$-NJI conditions. From the latter, the existence of a path of strict decay for every finite subnetwork follows. It remains open whether one can get back from the uniform NJI conditions to the existence of a path of strict decay for the original network. In the finite case, we obtain three different characterizations of the existence of a path of strict decay.%

In Section \ref{sec_classes}, we apply our results to four classes of gain operators.%

{\bf Subsection \ref{subsec_uniformly_continuous}}: We study gain operators which are uniformly continuous on their whole domain of definition and provide conditions under which UGAS implies the $\oplus$-MBI property, exploiting the structural similarity between the two uniform NJI conditions.%

{\bf Subsection \ref{subsec_maxtype}}: We study max-type gain operators, which have the property that they preserve the $\oplus$-operation on $\ell^{\infty}_+(\IC)$, i.e.~$\Gamma(s^1 \oplus s^2) = \Gamma(s^1) \oplus \Gamma(s^2)$. For these operators, we obtain a fairly complete theory. Using the results of Subsections \ref{subsec_necessary} and \ref{subsec_path_construction}, we can prove that the existence of a path of strict decay is equivalent to UGAS of $\Sigma(\Gamma_{\rho})$ for some $\rho \in \KC_{\infty}$, which in turn is equivalent to the combination of the uniform NJI conditions.%

{\bf Subsection \ref{subsec_homogeneous}}: We consider gain operators, which are homogeneous, i.e.~$\Gamma(\alpha s) = \alpha \Gamma(s)$ for $s \in \ell^{\infty}_+(\IC)$ and $\alpha \geq 0$, and subadditive, i.e.~$\Gamma(s^1 + s^2) \leq \Gamma(s^1) + \Gamma(s^2)$ for $s^1,s^2 \in \ell^{\infty}_+(\IC)$. For this class of operators, we show that the existence of a path of strict decay is equivalent to UGAS of $\Sigma(\Gamma)$, which in turn is equivalent to a (generalized) spectral radius condition. Moreover, we prove that under the restriction to linear gain functions for building enlarged gain operators also the $\oplus$-MBI property is an equivalent condition.%

{\bf Subsection \ref{subsec_finite_case}}: We consider the case when the index set $\IC$ is finite. Here, we recover and improve the known results from the literature. The improvement mainly consists in the replacement of the assumption of subadditivity for the functions defining the ISS formulation with uniform continuity.%

In Section \ref{sec_summary}, we present a list of open questions and problems for future work. Finally, the appendix, Section \ref{sec_appendix}, contains some useful technical lemmas.

\section{Gain operators}\label{sec_gain_ops}

In this section, we introduce the central objects of the paper and study their fundamental properties.%

{\bf Notation}: We write $\N = \{1,2,3,\ldots\}$ for the natural numbers, $\N_0 = \{0\} \cup \N$, and $\Z$ for the integers. If $a,b \in \Z$ with $a \leq b$, we write $[a;b] = \{ n \in \Z : a \leq n \leq b \}$. We further write $\R$ for the reals and $\R_+ = \{ r \in \R : r \geq 0 \}$. If $A$ is a finite set, $|A|$ denotes its cardinality. If $X,Y$ are topological spaces, we use the notation $C^0(X,Y)$ for the set of continuous mappings from $X$ into $Y$. For any $r \in \R$, we write $\lfloor r \rfloor = \max\{ n \in \Z : n \leq r \}$ and $\lceil r \rceil = \min\{ n \in \Z : r \leq n \}$.%

\subsection{Preliminaries}\label{subsec_prelim}

The spaces of $\KC$-functions and $\KC_{\infty}$-functions, respectively, are defined by%
\begin{align*}
  \KC &:= \left\{ \gamma \in C^0(\R_+,\R_+) : \gamma(0) = 0 \mbox{ and } \gamma \mbox{ is strictly increasing} \right\}, \\
  \KC_{\infty} &:= \left\{ \gamma \in \KC : \lim_{r \rightarrow \infty}\gamma(r) = \infty \right\}.%
\end{align*}
We further introduce the classes of $\LC$-functions and $\KC\LC$-functions:%
\begin{align*}
  \LC &:= \left\{ \gamma \in C^0(\R_+,\R_+) : \gamma \mbox{ is strictly decreasing with } \lim_{t\rightarrow\infty}\gamma(t) = 0 \right\}, \\
	\KC\LC &:= \left\{ \beta \in C^0(\R_+ \tm \R_+,\R_+) : \beta(\cdot,t) \in \KC,\ \forall t\geq0,\ \beta(r,\cdot) \in \LC,\ \forall r>0 \right\}.%
\end{align*}

For any nonempty countable index set $\IC$, we consider the ordered Banach space $(\ell^{\infty}(\IC),\ell^{\infty}_+(\IC))$, where%
\begin{align*}
  \ell^{\infty}(\IC) &:= \bigl\{ s = (s_i)_{i \in \IC} : s_i \in \R,\ \sup_{i\in\IC}|s_i| < \infty \bigr\}, \\
	\ell^{\infty}_+(\IC) &:= \left\{ s \in \ell^{\infty}(\IC) : s_i \geq 0,\ \forall i \in \IC \right\}.%
\end{align*}
The norm on $\ell^{\infty}(\IC)$ is the usual norm of uniform convergence:%
\begin{equation*}
  \|s\| := \sup_{i \in \IC}|s_i|.%
\end{equation*}
The cone $\ell^{\infty}_+(\IC)$ is closed, convex, and has nonempty interior given by%
\begin{equation*}
  \inner(\ell^{\infty}_+(\IC)) = \bigl\{ s \in \ell^{\infty}_+(\IC)\ :\ \exists \alpha > 0 \mbox{ with } \inf_{i\in\IC}s_i \geq \alpha \bigr\}.%
\end{equation*}
The partial order on $\ell^{\infty}(\IC)$, induced by the cone $\ell^{\infty}_+(\IC)$, is given by%
\begin{equation*}
  s^1 \leq s^2 \quad :\Leftrightarrow \quad s^2 - s^1 \in \ell^{\infty}_+(\IC).%
\end{equation*}
Hence, $s \in \ell^{\infty}_+(\IC)$ if and only if $s \geq 0$. We write $s^1 < s^2$ if $s^1 \leq s^2$ and $s^1 \neq s^2$. We further write $s^1 \ll s^2$ if $s^2 - s^1 \in \inner(\ell^{\infty}_+(\IC))$. For two subsets $A,B \subset \ell^{\infty}_+(\IC)$, we write $A \leq B$ if $a \leq b$ for all $(a,b) \in A \tm B$. For any $0 \leq s^1 \leq s^2$,  we define the \emph{order interval} $[s^1,s^2] := \{ s \in \ell^{\infty}_+(\IC) : s^1 \leq s \leq s^2 \}$.%

The Banach space $\ell^{\infty}(\IC)$ can be identified with the topological dual of $\ell^1(\IC)$ (the space of absolutely summable sequences). In this identification, any $s \in \ell^{\infty}(\IC)$ acts linearly on elements of $\ell^1(\IC)$ by%
\begin{equation}\label{eq_dual_action}
  s(t) := \sum_{i \in \IC} s_it_i \mbox{\quad for all\ } t = (t_i)_{i\in\IC} \in \ell^1(\IC).%
\end{equation}
On every topological dual space $X^*$ of a normed space $X$, we have the weak$^*$-topology which is the coarsest topology such that for every $x \in X$ the linear functional $x^* \mapsto x^*(x)$, defined on $X^*$, is continuous. If a sequence $(s^n)_{n\in\N}$ in $\ell^{\infty}(\IC)$ converges in the weak$^*$-topology to some $s$, we write $s^n \stackrel{\star}{\rightharpoonup} s$. By Lemma \ref{lem_weakstar1}, this is the case if and only if $(s^n)_{n\in\N}$ is norm-bounded and converges to $s$ componentwise. The weak$^*$-topology on $\ell^{\infty}(\IC)$ coincides with the norm-topology if and only if $\IC$ is finite.%
 
An operator $T:\ell^{\infty}_+(\IC) \rightarrow \ell^{\infty}_+(\IC)$ is called \emph{monotone} if $s^1 \leq s^2$ implies $T(s^1) \leq T(s^2)$ for all $s^1,s^2 \in \ell^{\infty}_+(\IC)$. The discrete-time dynamical system generated by $T$ is given by%
\begin{equation*}
  \Sigma(T):\quad s^{n+1} = T(s^n),\quad s^0 \in \ell^{\infty}_+(\IC),\ n \in \N_0.%
\end{equation*}
For any $s \in \ell^{\infty}_+(\IC)$, the sequence $(T^n(s))_{n\in\N_0}$ is called the \emph{trajectory of $\Sigma(T)$ starting in $s$}. Under the assumption that $T(0) = 0$, we call $\Sigma(T)$%
\begin{itemize}
\item \emph{uniformly globally stable (UGS)} if there exists $\varphi \in \KC_{\infty}$ such that $\|T^n(s)\| \leq \varphi(\|s\|)$ for all $s \in \ell^{\infty}_+(\IC)$ and $n \in \N_0$.%
\item \emph{uniformly globally asymptotically stable (UGAS)} if there exists $\beta \in \KC\LC$ such that $\|T^n(s)\| \leq \beta(\|s\|,n)$ for all $s \in \ell^{\infty}_+(\IC)$ and $n \in \N_0$.%
\item \emph{globally attractive (GATT)} if every trajectory of $\Sigma(T)$ converges to the origin in the norm-topology.%
\item \emph{uniformly globally attractive (UGATT)} if for all $r,\ep > 0$ there exists $n \in \N$ such that%
\begin{equation*}
  \|s\| \leq r \mbox{\ and\ } k \geq n \quad\Rightarrow\quad \|T^k(s)\| \leq \ep.%
\end{equation*}
\item \emph{weakly$^*$ globally attractive (GATT$^{\,*}$)} if every trajectory of $\Sigma(T)$ converges to the origin in the weak$^*$-topology.%
\end{itemize}
It is clear that UGAS $\Rightarrow$ UGATT $\Rightarrow$ GATT $\Rightarrow$ GATT$^*$ and UGAS $\Rightarrow$ UGS. From Lemma \ref{lem_ugas}, it follows that UGAS $\Leftrightarrow$ UGS $\wedge$ GATT.%

A point $s \in \ell^{\infty}_+(\IC)$ is called a \emph{point of decay} for $T$ if $T(s) \leq s$ and a \emph{fixed point} if $T(s) = s$. The \emph{set of decay} $\Psi(T)$ and the \emph{fixed point set} $\Fix(T)$ of $T$ are defined as the sets of all points of decay and fixed points, respectively.%

By $e^i$, $i \in \IC$, we denote the unit vectors in $\ell^{\infty}(\IC)$, i.e.~$e^i_j = \delta_{ij}$ with the Kronecker delta $\delta_{ij}$. By $\unit$, we denote the vector all of whose components are $1$s, i.e.~$\unit = \sum_{i\in\IC}e^i$.%

A natural binary operation on $\ell^{\infty}_+(\IC)$ is the componentwise maximum,%
\begin{equation*}
  s^1 \oplus s^2 := \left(\max\{s^1_i,s^2_i\}\right)_{i \in \IC} \mbox{\quad for all\ } s^1,s^2 \in \ell^{\infty}_+(\IC).%
\end{equation*}
It can be extended to arbitrary norm-bounded families of vectors $s^a$, $a \in A$, by%
\begin{equation*}
  \bigoplus_{a \in A}s^a := \Bigl(\sup_{a\in A}s^a_i\Bigr)_{i \in \IC}.%
\end{equation*}

A set $A \subset \ell^{\infty}_+(\IC)$ is called \emph{cofinal} if for every $s \in \ell^{\infty}_+(\IC)$ there is $\hat{s} \in A$ such that $s \leq \hat{s}$. It is called \emph{uniformly cofinal} if there exists $\varphi \in \KC_{\infty}$ such that for every $s \in \ell^{\infty}_+(\IC)$ there is $\hat{s} \in A$ with $s \leq \hat{s}$ and $\|\hat{s}\| \leq \varphi(\|s\|)$.%

A sequence $(s^n)_{n \in \N}$ in $\ell^{\infty}_+(\IC)$ is called \emph{increasing} if $s^n \leq s^{n+1}$ and \emph{decreasing} if $s^{n+1} \leq s^n$ for all $n \in \N$.%

We call a set $A \subset \ell^{\infty}_+(\IC)$ \emph{coercive} if there exists $\underline{\varphi} \in \KC_{\infty}$ such that $s \geq \underline{\varphi}(\|s\|)\unit$ for all $s \in A$. This inequality can also be written as%
\begin{equation*}
  \inf_{i\in\IC}s_i \geq \underline{\varphi}(\sup_{i\in\IC}s_i).%
\end{equation*}
Hence, a set $A$ is coercive if the variation in the components of each $s \in A$ is limited by a $\KC_{\infty}$-function.%

For any subset $\JC \subset \IC$, we define%
\begin{equation*}
  \ell^{\infty}_+(\IC,\JC) := \{ s \in \ell^{\infty}_+(\IC) : s_i = 0,\ \forall i \in \IC \backslash \JC \}%
\end{equation*}
which can be identified isometrically with $\ell^{\infty}_+(\JC)$. For every $s \in \ell^{\infty}_+(\IC)$, we let%
\begin{equation*}
  s_{|\JC} := \sum_{i\in \JC}s_ie^i \in \ell^{\infty}_+(\IC,\JC).%
\end{equation*}

\subsection{Definition and elementary properties}\label{subsec_def_props}

We now define gain operators as specific monotone operators on $\ell^{\infty}_+(\IC)$. A gain operator is associated with a network of dynamical systems, indexed by $\IC$, where each system is directly influenced by only finitely many neighbors. The set of indices of the neighbors of system $i$ is denoted by $\IC_i \subset \IC$. Here, we assume that $i \notin \IC_i$ and do \emph{not} exclude the case $\IC_i = \emptyset$. We also define an associated \emph{interconnection graph} $\GC$, which is a directed graph with vertex set $V(\GC)$ and edge set $E(\GC)$ given by%
\begin{equation*}
  V(\GC) := \IC \mbox{\quad and \quad} E(\GC) := \{ ji : i \in \IC, j \in \IC_i \},%
\end{equation*}
respectively. For a vertex $i$ and a positive integer $n \in \N$, we denote by $\NC^-_i(n)$ the set of all vertices $j$ such that there exists a directed path from $j$ to $i$ of length at most $n$. In contrast, the set $\NC^+_i(n)$ is the set of all vertices $j$ such that there exists a directed path from $i$ to $j$ of length at most $n$. In addition, we put $\NC^-_i(0) = \NC^+_i(0) := \{i\}$ and define%
\begin{equation*}
  \NC^-_i := \bigcup_{n \in \N_0}\NC^-_i(n),\quad \NC^+_i := \bigcup_{n \in \N_0}\NC^+_i(n).%
\end{equation*}
Observe that $\NC^-_i(1) = \{i\} \cup \IC_i$. Inductively, it follows that the sets $\NC^-_i(n)$ are all finite. In contrast, the sets $\NC^+_i(n)$ may be infinite.%

A gain operator compatible with the directed graph $\GC$ can be built from two families of functions. The first family consists of functions $\gamma_{ij} \in \KC_{\infty}$, $ji \in E(\GC)$, called \emph{interconnection gains}. We assume that this family is pointwise equicontinuous, i.e.~for all $r_0 \in \R_+$ and $\ep > 0$ there exists $\delta > 0$ such that $|r - r_0| \leq \delta$ implies $|\gamma_{ij}(r) - \gamma_{ij}(r_0)| \leq \ep$ for all $ji \in E(\GC)$ and $r \in \R_+$.%

The second family has elements $\mu_i:\ell^{\infty}_+(\IC) \rightarrow [0,\infty]$, $i \in \IC$, called \emph{monotone aggregation functions (MAFs)},\footnote{This term is adopted from the small-gain literature.} satisfying the following properties:%
\begin{enumerate}
\item[(M1)] \emph{Uniform positive definiteness}: There exists $\xi \in \KC_{\infty}$ such that $\mu_i(0) = 0$ and $\mu_i(s) \geq \xi(\|s\|)$ for all $s \in \ell^{\infty}_+(\IC)$ and $i \in \IC$.%
\item[(M2)] \emph{Monotonicity}: $0 \leq s^1 \leq s^2$ implies $\mu_i(s^1) \leq \mu_i(s^2)$ for all $i \in \IC$.%
\item[(M3)] \emph{Finite continuity}: For any finite set $\JC \subset \IC$ and $i \in \IC$, the restriction of $\mu_i$ to $\ell^{\infty}_+(\IC,\JC)$ is finite-valued and continuous.%
\item[(M4)] \emph{Equicontinuity}: For every norm-bounded set $A \subset \ell^{\infty}_+(\IC)$ and $\ep > 0$, there is $\delta > 0$ such that for all $s^0 \in A$ and $s \in \ell^{\infty}_+(\IC)$:%
\begin{equation*}
  \|s - s^0\| \leq \delta \quad \Rightarrow \quad \sup_{i \in \IC}|\mu_i(s_{|\IC_i}) - \mu_i(s^0_{|\IC_i})| \leq \ep.%
\end{equation*}
\end{enumerate}

We now define gain operators.%

\begin{definition}\label{def_normal_gain_op}
The \emph{gain operator} $\Gamma:\ell^{\infty}_+(\IC) \rightarrow \ell^{\infty}_+(\IC)$ associated with the graph $\GC$, the interconnection gains $\gamma_{ij}$, and the MAFs $\mu_i$, is defined by%
\begin{equation*}
  \Gamma(s) = (\Gamma_i(s))_{i \in \IC},\quad \Gamma_i(s) := \mu_i([\gamma_{ij}(s_j)]_{j\in\IC_i}),%
\end{equation*}
where we use the natural identification of $\ell^{\infty}_+(\IC_i)$ with $\ell^{\infty}_+(\IC,\IC_i)$.
\end{definition}

In the rest of the paper, we speak of the \emph{finite case} when we consider gain operators with $|\IC| < \infty$ and otherwise of the \emph{infinite case}.%

\begin{example}
We consider two of the most relevant examples for MAFs. Let%
\begin{equation*}
  \mu_{\oplus}(s) := \sup_{i\in\IC}s_i = \|s\|,\quad \mu_{\oplus}:\ell^{\infty}_+(\IC) \rightarrow [0,\infty].%
\end{equation*}
In case that $\mu_i \equiv \mu_{\oplus}$, we call $\Gamma$ a \emph{max-type operator}. Here, (M1) holds with $\xi = \id$ and (M2), (M3) are trivial. If $\|s - s^0\| \leq \delta$ for some $s^0,s \in \ell^{\infty}_+(\IC)$, the triangle inequality from below implies%
\begin{equation*}
  \sup_{i \in \IC} \Bigl|\sup_{j \in \IC_i}s_j - \sup_{j \in \IC_i}s^0_j\Bigr| \leq \sup_{ji \in E(\GC)} |s_j - s^0_j| \leq \|s - s^0\| \leq \delta.%
\end{equation*}
Hence, also (M4) is satisfied (even on the whole cone $\ell^{\infty}_+(\IC)$). We thus see that for max-type operators the only requirement is that the interconnection gains form a pointwise equicontinuous family. Now, consider%
\begin{equation*}
  \mu_{\otimes}(s) := \sum_{i\in\IC}s_i,\quad \mu_{\otimes}:\ell^{\infty}_+(\IC) \rightarrow [0,\infty].%
\end{equation*}
In case that $\mu_i \equiv \mu_{\otimes}$, we call $\Gamma$ a \emph{sum-type operator}. Here, (M1) again holds with $\xi = \id$ and (M2), (M3) are trivial. Concerning (M4), $\|s - s^0\| \leq \delta$ implies%
\begin{equation*}
  \Bigl|\sum_{j \in \IC_i}s_j - \sum_{j \in \IC_i}s^0_j\Bigr| \leq \sum_{j \in \IC_i} |s_j - s^0_j| \leq \delta \cdot |\IC_i|,%
\end{equation*}
which can hold with equality in certain cases. Hence, for (M4) it is necessary and sufficient that the cardinality $|\IC_i|$ is bounded over $i \in \IC$.
\end{example}

Elementary properties of gain operators are summarized in the following result.%

\begin{proposition}\label{prop_gainop_props}
Any gain operator $\Gamma$ has the following properties:%
\begin{enumerate}
\item[(a)] $\Gamma$ is well-defined, i.e.~$\Gamma(s) \in \ell^{\infty}_+(\IC)$ for all $s \in \ell^{\infty}_+(\IC)$.%
\item[(b)] $\Gamma(0) = 0$ and $\Gamma$ is a monotone operator.
\item[(c)] $\Gamma$ is continuous and uniformly continuous on bounded sets.%
\item[(d)] $\Gamma$ is sequentially continuous in the weak$^*$-topology.%
\end{enumerate}
\end{proposition}

\begin{proof}
First, consider the family $\FC := \{ f_i : i \in \IC \} \subset C^0(\R_+,\R_+)$ with $f_i(r) :\equiv \mu_i(r\unit_{|\IC_i})$. By (M3), each $f_i$ is finite-valued and continuous. We claim that $\FC$ is uniformly bounded on every compact interval $[0,r_0]$. To show this, we use (M4) with $\ep = 1$, which gives $\delta > 0$ such that $|r - r'| \leq \delta$ with $r,r' \in [0,r_0]$ implies%
\begin{equation*}
  \sup_{i \in \IC}|f_i(r') - f_i(r)| = \sup_{i \in \IC}\left|\mu_i\left(r'\unit_{|\IC_i}\right) - \mu_i\left(r\unit_{|\IC_i}\right)\right| \leq 1,%
\end{equation*}
leading to $f_i(r_0) = |f_i(r_0) - f_i(0)| \leq \lceil r_0/\delta \rceil$ for all $i \in \IC$.%

Now, we show that $\Gamma$ is well-defined. To this end, let $s \in \ell^{\infty}_+(\IC)$ and put $r := \|s\|$. Then, by the monotonicity of each $\mu_i$ and $\gamma_{ij}$, we have%
\begin{equation*}
  \Gamma_i(s) = \mu_i( [\gamma_{ij}(s_j)]_{j\in\IC_i} ) \leq \mu_i( [\gamma_{ij}(r)]_{j\in\IC_i} ) \mbox{\quad for all\ } i \in \IC.%
\end{equation*}
With the same arguments as used above for the family $\FC$, the equicontinuity of $\{\gamma_{ij}\}$ implies the existence of $R > 0$ with $\gamma_{ij}(r) \leq R$ for all $ji \in E(\GC)$. We conclude the existence of $B > 0$ with%
\begin{equation*}
  \|\Gamma(s)\| = \sup_{i\in\IC}\Gamma_i(s) \leq \sup_{i\in\IC} \mu_i\left( R\unit_{|\IC_i} \right) = \sup_{i\in\IC} f_i(R) \leq B,%
\end{equation*}
showing that $\Gamma$ is well-defined. From the monotonicity assumptions on $\gamma_{ij}$ and $\mu_i$, it follows that $\Gamma$ is a monotone operator. By (M1), it is clear that $\Gamma(0) = 0$. Hence, (a) and (b) are proven.%

We now prove the uniform continuity of $\Gamma$ on bounded sets. To this end, let $B \subset \ell^{\infty}_+(\IC)$ be bounded, say $\|s\| \leq r$ for all $s \in B$. Then also the set%
\begin{equation*}
  A := \left\{ (\gamma_{ij}(s_j))_{j \in \IC_i} \in \ell^{\infty}_+(\IC,\IC_i)\ :\ s \in B,\ i \in \IC \right\}%
\end{equation*}
is bounded, because the equicontinuity of $\{\gamma_{ij}\}$ implies $\gamma_{ij}(s_j) \leq \gamma_{ij}(r) \leq R$ for all $ji \in E(\GC)$ and some $R > 0$. For a given $\ep > 0$, we first choose $\delta = \delta(A,\ep)$ according to (M4). We subsequently choose $\rho > 0$ so that $\|s^1 - s^2\| \leq \rho$ with $s^1,s^2 \in B$ implies $|\gamma_{ij}(s^1_j) - \gamma_{ij}(s^2_j)| \leq \delta$ for all $ji \in E(\GC)$, which is possible by the pointwise equicontinuity of $\{\gamma_{ij}\}$ (which is uniform on compact intervals). Altogether, $\|s^1 - s^2\| \leq \rho$ with $s^1,s^2 \in B$ implies%
\begin{equation*}
  \|\Gamma(s^1) - \Gamma(s^2)\| = \sup_{i \in \IC} \left|\mu_i( [\gamma_{ij}(s^1_j)]_{j \in \IC_i} ) - \mu_i( [\gamma_{ij}(s^2_j)]_{j \in \IC_i} )\right| \leq \ep.%
\end{equation*}
Hence, (c) is proved. To prove (d), consider a sequence $(s^n)_{n\in\N}$ with $s^n \stackrel{\star}{\rightharpoonup} s$ for some $s$. By Lemma \ref{lem_weakstar1}, this means that $(s^n)$ is norm-bounded and converges componentwise to $s$. Then, for every $i \in \IC$ we have%
\begin{align*}
  \lim_{n \rightarrow \infty} \Gamma_i(s^n) &= \lim_{n \rightarrow \infty} \mu_i( [\gamma_{ij}(s^n_j)]_{j\in\IC_i} ) = \mu_i\left( \lim_{n \rightarrow \infty} [\gamma_{ij}(s^n_j)]_{j\in\IC_i} \right) \\
																						&= \mu_i\Bigl( \Bigl[\lim_{n\rightarrow\infty}\gamma_{ij}(s^n_j)\Bigr]_{j\in\IC_i} \Bigr) = \mu_i\Bigl( \Bigl[\gamma_{ij}\Bigl(\lim_{n\rightarrow\infty}s^n_j\Bigr)\Bigr]_{j\in\IC_i} \Bigr) \\
																						&= \mu_i\left( \left[\gamma_{ij}(s_j)\right]_{j\in\IC_i} \right) = \Gamma_i(s).%
\end{align*}
Here, the second equality follows from (M3). The third equality follows from finiteness of $\IC_i$, and the fourth from continuity of $\gamma_{ij}$. Hence, $(\Gamma(s^n))_{n\in\N}$ converges to $\Gamma(s)$ componentwise. Moreover, if $B > 0$ is a bound on $\|s^n\|$, then $s^n \leq B\unit$ and thus $\Gamma(s^n) \leq \Gamma(B\unit)$ due to monotonicity, implying $\|\Gamma(s^n)\| \leq \|\Gamma(B\unit)\|$ for all $n$. Consequently, $\Gamma(s^n) \stackrel{\star}{\rightharpoonup} \Gamma(s)$.
\end{proof}

\begin{remark}
The sequential continuity of $\Gamma$ in the weak$^*$-topology has the following important implication: If a norm-bounded trajectory $(\Gamma^n(s))_{n\in\N}$ converges componentwise to $s^* \in \ell^{\infty}_+(\IC)$, then $s^*$ is a fixed point of $\Gamma$. This will be used many times throughout the paper.
\end{remark}

Almost all results of this paper hold true for a more general class of operators that, for simplicity, we will also call gain operators.%

\begin{definition}
Given a directed graph $\GC$ as described above, an operator $\Gamma:\ell^{\infty}_+(\IC) \rightarrow \ell^{\infty}_+(\IC)$ is called a \emph{generalized gain operator (compatible with $\GC$)} if it satisfies the following properties:%
\begin{enumerate}
\item[(G1)] For each $s \in \ell^{\infty}_+(\IC)$ and $i \in \IC$, $\Gamma_i(s)$ only depends on the components $s_j$ with $j \in \IC_i$.%
\item[(G2)] $\Gamma(0) = 0$ and $\Gamma$ is a monotone operator.%
\item[(G3)] $\Gamma$ is uniformly continuous on bounded sets.%
\item[(G4)] $\Gamma$ is sequentially continuous in the weak$^*$-topology.%
\end{enumerate}
For generalized gain operators, we also allow $\IC_i$ to be an infinite set, possibly with $i \in \IC_i$. If all $\IC_i$ are finite, we speak of a \emph{locally finite generalized gain operator}.
\end{definition}

\begin{remark}
Observe that in the finite case a generalized gain operator is simply a continuous monotone operator that fixes the origin. Property (G3) holds, because the closure of a bounded set is compact here, and (G4) holds, because the weak$^*$-topology coincides with the norm-topology.
\end{remark}

Any function $\varphi \in \KC_{\infty}$ acts as a continuous monotone operator on $\ell^{\infty}_+(\IC)$ by%
\begin{equation*}
  \varphi(s) := (\varphi(s_i))_{i\in\IC},\quad \varphi:\ell^{\infty}_+(\IC) \rightarrow \ell^{\infty}_+(\IC).%
\end{equation*}
This operator has the following properties:%
\begin{itemize}
\item It is invertible with inverse $\varphi^{-1}(s) = (\varphi^{-1}(s_i))_{i\in\IC}$.%
\item $\|\varphi(s)\| = \varphi(\|s\|)$ for all $s \in \ell^{\infty}_+(\IC)$.%
\item $\varphi(s^1 \oplus s^2) = \varphi(s^1) \oplus \varphi(s^2)$ for all $s^1,s^2 \in \ell^{\infty}_+(\IC)$.%
\end{itemize}
Given a generalized gain operator $\Gamma$ and $\rho \in \KC_{\infty}$, we define two \emph{enlarged operators}%
\begin{equation*}
  \Gamma_{\rho} := (\id + \rho) \circ \Gamma \mbox{\quad and \quad} \Gamma^{\rho} := \Gamma \circ (\id + \rho).%
\end{equation*}
It is not hard to see that $\Gamma_{\rho}$ and $\Gamma^{\rho}$ are generalized gain operators as well. They are conjugate, since%
\begin{equation}\label{eq_op_conj}
  (\id + \rho) \circ \Gamma^{\rho} = \Gamma_{\rho} \circ (\id + \rho).%
\end{equation}
If $\Gamma$ is a gain operator, then $\Gamma_{\rho}$ is the gain operator associated with the gains $\{\gamma_{ij}\}$ and MAFs $\{(\id + \rho) \circ \mu_i\}$,\footnote{To prove that this family satisfies (M4), one can use that the $\mu_i$ are uniformly continuous on bounded sets and $\id + \rho$ is uniformly continuous on compact intervals.} while $\Gamma^{\rho}$ is the gain operator associated with the gains $\{\gamma_{ij} \circ (\id + \rho)\}$ and MAFs $\{\mu_i\}$.

A mapping $\sigma:\R_+ \rightarrow \ell^{\infty}_+(\IC)$ is called a \emph{path} in the cone $\ell^{\infty}_+(\IC)$. A path $\sigma$ is called \emph{increasing} if $r_1 \leq r_2$ implies $\sigma(r_1) \leq \sigma(r_2)$, and \emph{strictly increasing} if $r_1 < r_2$ implies $\sigma(r_1) \ll \sigma(r_2)$ for all $r_1,r_2 \in \R_+$.%

We now introduce the notion of a \emph{path of strict decay} (which is adopted from \cite[Def.~II.10]{KMZ}). The existence of a path of strict decay for a gain operator $\Gamma$, together with other assumptions, guarantees that an ISS Lyapunov function for the overall network can be constructed from ISS Lyapunov functions of the subsystems, which is described in detail in \cite[Sec.~2]{KZa}.%

\begin{definition}\label{def_path_of_decay}
Given a generalized gain operator $\Gamma$, a path $\sigma$ in the cone $\ell^{\infty}_+(\IC)$ is called a \emph{path of strict decay} for $\Gamma$ if it satisfies the following properties:%
\begin{enumerate}
\item[(i)] There exists $\rho \in \KC_{\infty}$ such that $\Gamma_{\rho}(\sigma(r)) \leq \sigma(r)$ for all $r \geq 0$.%
\item[(ii)] There exist $\varphi_{\min},\varphi_{\max} \in \KC_{\infty}$ such that $\varphi_{\min}(r)\unit \leq \sigma(r) \leq \varphi_{\max}(r)\unit$ for all $r \geq 0$.%
\item[(iii)] Each component function $\sigma_i$, $i \in \IC$, is a $\KC_{\infty}$-function.%
\item[(iv)] For every compact interval $K \subset (0,\infty)$, there exist constants $0 < l \leq L$ such that for all $r_1,r_2 \in K$ and $i \in \IC$%
\begin{equation*}
  l|r_1 - r_2| \leq |\sigma_i^{-1}(r_1) - \sigma_i^{-1}(r_2)| \leq L|r_1 - r_2|.%
\end{equation*}
\end{enumerate}
If all of the above properties hold with the exception that $\rho = 0$ in (i), we call $\sigma$ a \emph{path of decay} for $\Gamma$.
\end{definition}

\begin{remark}
Observe that a path of strict decay for $\Gamma$ is a path of decay for $\Gamma_{\rho}$ with $\rho$ as in (i). Using Proposition \ref{prop_path_of_decay}(a) below, it is easy to see that any path of decay is continuous, locally Lipschitz continuous on $(0,\infty)$, strictly increasing, and satisfies $\sigma(0) = 0$.
\end{remark}

The definition can be slightly simplified, using the next result.%

\begin{proposition}\label{prop_path_of_decay}
Regarding Definition \ref{def_path_of_decay}, the following holds:%
\begin{enumerate}
\item[(a)] In property (iv), the functions $\sigma_i^{-1}$ can be replaced by $\sigma_i$.%
\item[(b)] In property (ii), we can assume without loss of generality that $\varphi_{\min} = \id$.%
\end{enumerate}
\end{proposition}

\begin{proof}
(a) Assume that for every compact interval $L \subset (0,\infty)$ there are constants $0 < c \leq C$ such that%
\begin{equation*}
  c|r_1 - r_2| \leq |\sigma_i(r_1) - \sigma_i(r_2)| \leq C|r_1 - r_2|%
\end{equation*}
for all $i \in \IC$ and $r_1,r_2 \in L$. If $K = [a,b] \subset (0,\infty)$ is any compact interval, then%
\begin{equation*}
  \sigma_i^{-1}([a,b]) = [\sigma_i^{-1}(a),\sigma_i^{-1}(b)] \subset [\varphi_{\max}^{-1}(a),\varphi_{\min}^{-1}(b)] =: L \subset (0,\infty).%
\end{equation*}
Hence, with the corresponding constants $0 < c \leq C$, we obtain%
\begin{equation*}
  \frac{1}{C}|r_1 - r_2| \leq |\sigma_i^{-1}(r_1) - \sigma_i^{-1}(r_2)| \leq \frac{1}{c}|r_1 - r_2| \mbox{\quad for all\ } r_1,r_2 \in K.%
\end{equation*}

(b) Given a path of decay $\sigma$, we can choose $\tilde{\varphi}_{\min} \in \KC_{\infty}$ such that $\tilde{\varphi}_{\min} \leq \varphi_{\min}$ and both $\tilde{\varphi}_{\min}$, $\tilde{\varphi}_{\min}^{-1}$ are locally Lipschitz continuous on $(0,\infty)$; for instance, we can choose $z_1$ from Lemma \ref{lem_mir_A16} with $z_1 \leq \varphi_{\min}$ and put $\tilde{\varphi}_{\min}(r) :\equiv z_1(r)(1 - \rme^{-r})$. This $\KC_{\infty}$-function is smooth on $(0,\infty)$ and its derivative vanishes nowhere, implying that also its inverse is smooth on $(0,\infty)$. Then, we consider the reparametrized path $\tilde{\sigma} := \sigma \circ \tilde{\varphi}_{\min}^{-1}$. Clearly, $\tilde{\sigma}$ satisfies properties (i) and (iii). Moreover, by construction, $r\unit \leq \tilde{\sigma}(r) \leq \varphi_{\max} \circ \tilde{\varphi}_{\min}^{-1}(r)\unit$ for all $r \geq 0$. By the choice of the reparametrization, $\tilde{\sigma}$ also satisfies property (iv), and thus is a path of strict decay for $\Gamma$.
\end{proof}

The following proposition is the main result of this subsection. With this result at hand, finding a path of strict decay becomes a considerably easier task.%

\begin{proposition}\label{prop_path_of_decay_2}
Given a generalized gain operator $\Gamma$, assume that there exists a continuous and increasing path $\sigma$ in $\ell^{\infty}_+(\IC)$, satisfying properties (i) and (ii) in Definition \ref{def_path_of_decay}. Then there exists a path of strict decay for $\Gamma$.
\end{proposition}

\begin{proof}
The proof consists of two steps. In the first step, we assume that $\sigma$ is strictly increasing, while in the second step we weaken this assumption.%

\emph{Step 1}: Let $\sigma$ be a path satisfying (i) and (ii) in Definition \ref{def_path_of_decay}, which is continuous and strictly increasing. Consider a compact interval $[a,b] \subset (0,\infty)$. We show that the restriction of $\sigma$ to $[a,b]$ can be approximated by a piecewise linear path evolving in $\Psi(\Gamma)$. The fact that such an approximation is possible will yield a piecewise linear path of decay. By Lemma \ref{lem_kinfty}, there exists $\eta \in \KC_{\infty}$ with $(\id + \rho)^{-1} = \id - \eta$. By assumption, $\sigma$ is continuous, hence uniformly continuous on the compact interval $[a,b]$. Thus, we can find $\delta > 0$ such that $|r_1 - r_2| < \delta$ implies $\|\sigma(r_1) - \sigma(r_2)\| < \ep := \eta(\varphi_{\min}(a))$, whenever $r_1,r_2 \in [a,b]$. Now, let $s^{\alpha} := (1 - \alpha)\sigma(r_1) + \alpha\sigma(r_2)$ for $\alpha \in [0,1]$ and $r_1,r_2 \in [a,b]$ with $r_1 < r_2 < r_1 + \delta$. This implies $\sigma(r_1) \leq s^{\alpha} \leq \sigma(r_2)$, and thus%
\begin{align*}
  \Gamma(s^{\alpha}) &\leq \Gamma(\sigma(r_2)) \leq (\id + \rho)^{-1}(\sigma(r_2)) = \sigma(r_2) - \eta(\sigma(r_2)) \\
	                   &\leq \sigma(r_2) - \eta(\sigma(a)) \leq \sigma(r_2) - \eta(\varphi_{\min}(a))\unit = \sigma(r_2) - \ep\unit \leq \sigma(r_1) \leq s^{\alpha}.%
\end{align*}
From the assumption that $r_1 < r_2$ implies $\sigma(r_1) \ll \sigma(r_2)$, we conclude that the slope of $\alpha \mapsto s^{\alpha}_i$ has a uniform lower bound. In fact, if $\sigma(r_2) - \sigma(r_1) \geq \gamma\unit$, then $(\partial s^{\alpha}_i)/(\partial \alpha) = \sigma_i(r_2) - \sigma_i(r_1) \geq \gamma$ for all $i \in \IC$. Moreover, $(\partial s_i^{\alpha})/(\partial \alpha) \leq \|\sigma(r_2) - \sigma(r_1)\|$. By exhausting $(0,\infty)$ with compact intervals, concatenating the piecewise linear paths obtained in this way, and properly reparametrizing them, we obtain a continuous path $\tilde{\sigma}:\R_+ \rightarrow \ell^{\infty}_+(\IC)$, which is piecewise linear on $(0,\infty)$, evolves in $\Psi(\Gamma)$, and satisfies (ii), (iii) and (iv). More precisely, we can find a strictly increasing double-sided sequence $(r_n)_{n\in\Z}$ in $(0,\infty)$ such that $r_n \rightarrow 0$ for $n \rightarrow -\infty$ and $r_n \rightarrow \infty$ for $n \rightarrow \infty$, pick for each $r > 0$ the unique $n = n_r \in \Z$ with $r \in [r_n,r_{n+1})$, and define $\tilde{\sigma}$ by%
\begin{equation*}
  \tilde{\sigma}(r) := \left(1 - \frac{r - r_n}{r_{n+1} - r_n}\right)\sigma(r_n) + \frac{r - r_n}{r_{n+1} - r_n}\sigma(r_{n+1}) \mbox{\quad for all\ } r > 0,%
\end{equation*}
extended by $\tilde{\sigma}(0) := 0$. Observe that upper and lower $\KC_{\infty}$-bounds for the components of $\tilde{\sigma}$ are obtained via the analogous piecewise linear approximation with $\sigma(r_n)$ replaced by $\varphi_{\min}(r_n)$ or $\varphi_{\max}(r_n)$, respectively. In this way, we obtain a path of decay for $\Gamma$. If we write $\id + \rho = (\id + \rho_2) \circ (\id + \rho_1)$ with $\rho_1,\rho_2 \in \KC_{\infty}$, using Lemma \ref{lem_kinfty}, we can replace $\Gamma$ with $\Gamma_{\rho_1}$ in our proof, and thus obtain a path of decay for $\Gamma_{\rho_1}$, which is a path of strict decay for $\Gamma$.%

\emph{Step 2}: Now, let $\sigma$ be a continuous and increasing path satisfying (i) and (ii) in Definition \ref{def_path_of_decay}. First, observe that the path $\tilde{\sigma} := (\id + \rho)^{-1} \circ \sigma$ satisfies%
\begin{equation*}
  \Gamma^{\rho} \circ \tilde{\sigma} = \Gamma^{\rho} \circ (\id + \rho)^{-1} \circ \sigma \stackrel{\eqref{eq_op_conj}}{=} (\id + \rho)^{-1} \circ \Gamma_{\rho} \circ \sigma \leq (\id + \rho)^{-1} \circ \sigma = \tilde{\sigma}.%
\end{equation*}
Moreover, $\tilde{\sigma}$ clearly satisfies all other assumptions of this proposition. Hence, we will assume that $\sigma$ satisfies $\Gamma^{\rho} \circ \sigma \leq \sigma$ instead of $\Gamma_{\rho} \circ \sigma \leq \sigma$. By Step 1, it now suffices to show that there exists a continuous and strictly increasing path $\tilde{\sigma}$ satisfying (i) and (ii) in Definition \ref{def_path_of_decay}. This is done in four substeps.%

\emph{Step 2.a}: For each $r \geq 0$, we factorize $\id + \rho = (\id + \rho_1^r) \circ (\id + \rho_2^r)$ with $\rho_1^r,\rho_2^r \in \KC_{\infty}$, such that the following holds:%
\begin{enumerate}
\item[(a)] If $r_1 < r_2$, then $\rho_2^{r_1}(r) < \rho_2^{r_2}(r)$ for every $r > 0$.%
\item[(b)] The function $\rho_2:(r_1,r_2) \mapsto \rho_2^{r_1}(r_2)$ is continuous.%
\item[(c)] There exists $\tilde{\rho} \in \KC_{\infty}$ such that $\tilde{\rho} \leq \rho_2^r$ for all $r \geq 0$.%
\end{enumerate}
To achieve this, first choose a continuous and strictly increasing function $f:\R_+ \rightarrow (0,1)$. Then define%
\begin{equation*}
  \rho_2^r:= f(r)\rho \mbox{\quad for all\ } r \in \R_+.%
\end{equation*}
Clearly, $\rho_2^r \in \KC_{\infty}$ and $\rho_2^r(\tilde{r}) < \rho(\tilde{r})$ for all $r \geq 0$ and $\tilde{r} > 0$. We put%
\begin{equation*}
  \rho_1^r := (1 - f(r))\rho \circ (\id + \rho_2^r)^{-1}.%
\end{equation*}
Then clearly $\rho_1^r \in \KC_{\infty}$ and%
\begin{equation*}
  (\id + \rho_1^r) \circ (\id + \rho_2^r) = \id + f(r)\rho + (1 - f(r)) \rho = \id + \rho.%
\end{equation*}
Item (a) above holds, because $f$ is strictly increasing. The continuity of $\rho_2$ as a function of both arguments easily follows from the continuity of $\rho$ and $f$, i.e.~(b) holds. If we choose $f$ such that $f(r) \geq \alpha > 0$ for some $\alpha > 0$ (which is easily possible), it follows that $\rho_2^r \geq \alpha\rho =: \tilde{\rho}$ for all $r$, i.e.~(c) holds.%

\emph{Step 2.b}: We define the path $\tilde{\sigma}(r) := (\id + \rho_2^r)(\sigma(r))$ for all $r \geq 0$. Then%
\begin{align*}
  \Gamma_{\tilde{\rho}}(\tilde{\sigma}(r)) &= (\id + \tilde{\rho}) \circ \Gamma \circ (\id + \rho_2^r) \circ \sigma(r) \\
	                                         &\leq (\id + \rho_2^r) \circ \Gamma \circ (\id + \rho) \circ \sigma(r) = (\id + \rho_2^r) \circ \Gamma^{\rho} \circ \sigma(r) \\
																					 &\leq (\id + \rho_2^r) \circ \sigma(r) = \tilde{\sigma}(r).%
\end{align*}
Hence, $\tilde{\sigma}$ satisfies property (i) of a path of strict decay. Moreover, it holds that%
\begin{equation*}
  \varphi_{\min}(r)\unit \leq \sigma(r) \leq \tilde{\sigma}(r) \leq (\id + \rho) \circ \sigma(r) \leq (\id + \rho) \circ \varphi_{\max}(r)\unit,%
\end{equation*}
showing that $\tilde{\sigma}$ also satisfies property (ii) of a path of strict decay.%

\emph{Step 2.c}: Let $r_1 < r_2$ and assume without loss of generality that $r_1 > 0$. We then have $\sigma(r_1) \geq \varphi_{\min}(r_1)\unit \gg 0$. In particular, all components of the vector $\sigma(r_1)$ are contained in the compact interval $J := [\varphi_{\min}(r_1),\varphi_{\max}(r_1)] \subset (0,\infty)$. Observe that there exists $\alpha > 0$ with $\rho_2^{r_2}(r) - \rho_2^{r_1}(r) = (f(r_2) - f(r_1))\rho(r) \geq \alpha$ for all $r \in J$, implying that $\tilde{\sigma}$ is strictly increasing:%
\begin{equation*}
  \tilde{\sigma}(r_1) = (\id + \rho_2^{r_1})(\sigma(r_1)) \ll (\id + \rho_2^{r_2})(\sigma(r_1)) \leq (\id + \rho_2^{r_2})(\sigma(r_2)) = \tilde{\sigma}(r_2).%
\end{equation*}

\emph{Step 2.d:} It remains to prove the continuity of $\tilde{\sigma}$. To this end, fix $r_0 \geq 0$ and $\ep > 0$. We first observe that for any $r \geq 0$ it holds that%
\begin{equation*}
  \|\tilde{\sigma}(r) - \tilde{\sigma}(r_0)\| \leq \|\sigma(r) - \sigma(r_0)\| + \sup_{i \in \IC}|\rho_2(r,\sigma_i(r)) - \rho_2(r_0,\sigma_i(r_0))|.%
\end{equation*}
By continuity of $\sigma$, we can find $\delta > 0$ such that $|r - r_0| \leq \delta$ implies $\|\sigma(r) - \sigma(r_0)\| \leq \ep/2$. By uniform continuity of $\rho_2$ on compact subsets of $\R_+ \tm \R_+$, we see that the second term can also be made smaller than $\ep/2$ if $\delta$ is chosen small enough. Hence, $\tilde{\sigma}$ is continuous.
\end{proof}

Motivated by the preceding proposition, we introduce the following definition.%

\begin{definition}\label{def_c0_path}
Given a generalized gain operator $\Gamma$, a path $\sigma$ in the cone $\ell^{\infty}_+(\IC)$ is called a \emph{$C^0$-path of strict decay} for $\Gamma$ if it satisfies the following properties:%
\begin{enumerate}
\item[(i)] \emph{Strict decay}: There is $\rho \in \KC_{\infty}$ such that $\sigma(r) \in \Psi(\Gamma_{\rho})$ for all $r \geq 0$.%
\item[(ii)] \emph{Coercivity}: There are $\varphi_{\min},\varphi_{\max} \in \KC_{\infty}$ such that $\varphi_{\min}(r)\unit \leq \sigma(r) \leq \varphi_{\max}(r)\unit$ for all $r \geq 0$.%
\item[(iii)] \emph{Monotonicity}: $\sigma$ is increasing.%
\item[(iv)] \emph{Continuity}: $\sigma$ is continuous.
\end{enumerate}
If all of the above properties hold with the exception that $\rho = 0$ in (i), we call $\sigma$ a \emph{$C^0$-path of decay} for $\Gamma$.
\end{definition}

\subsection{Associated operators}\label{subsec_associated_ops}

Given a generalized gain operator $\Gamma$, we introduce further monotone operators associated with $\Gamma$, which are useful for analyzing under which conditions a $C^0$-path of strict decay exists. First, we define the \emph{augmented gain operator}%
\begin{equation*}
  \hat{\Gamma}(s) := s \oplus \Gamma(s),\quad \hat{\Gamma}:\ell^{\infty}_+(\IC) \rightarrow \ell^{\infty}_+(\IC).%
\end{equation*}
Moreover, for each $b \in \ell^{\infty}_+(\IC)$, we define the \emph{projected gain operator}\footnote{Observe that the mapping $s \mapsto b \oplus s$ can be regarded as a projection to the set $\{ s \geq b \}$.}%
\begin{equation*}
	\bsign{\Gamma}_b(s) := b \oplus \Gamma(s),\quad \bsign{\Gamma}_b:\ell^{\infty}_+(\IC) \rightarrow \ell^{\infty}_+(\IC).%
\end{equation*}

The properties of these operators listed in the following proposition can easily be derived, see \cite[Prop.~3.1 and Lem.~3.2]{KZa}.%

\begin{proposition}\label{prop_hatgamma_properties}
The following statements hold for all $b,s \in \ell^{\infty}_+(\IC)$:%
\begin{enumerate}
\item[(a)] $\bsign{\Gamma}_b$ and $\hat{\Gamma}$ are well-defined norm-continuous monotone operators.%
\item[(b)] $\bsign{\Gamma}_b$ and $\hat{\Gamma}$ are sequentially continuous in the weak$^*$-topology.%
\item[(c)] $\bsign{\Gamma}_b(0) = b$ and $\hat{\Gamma}(0) = 0$.%
\item[(d)] $\bsign{\Gamma}_b(s) \geq b$ and $\hat{\Gamma}(s) \geq s$.%
\item[(e)] $\Psi(\bsign{\Gamma}_b) \subset \Psi(\Gamma)$ and $\Fix(\hat{\Gamma}) = \Psi(\Gamma)$.%
\item[(f)] $\hat{\Gamma}^n(s) = \bsign{\Gamma}_s^n(s)$ for all $n \in \N_0$.%
\end{enumerate}
\end{proposition}

Observe that, due to (a) and (d) in Proposition \ref{prop_hatgamma_properties}, every trajectory of $\Sigma(\hat{\Gamma})$ is increasing. Hence, if such a trajectory is norm-bounded, due to (b) and (e) it weakly$^*$ converges to a point in $\Psi(\Gamma)$. To effectively bound all trajectories of $\Sigma(\hat{\Gamma})$, we simply assume that $\Sigma(\hat{\Gamma})$ is UGS. As the next proposition shows, this assumption is also necessary for the existence of a path of decay.%

\begin{proposition}\label{prop_hatsigma_ugs}
The following statements are equivalent:%
\begin{enumerate}
\item[(a)] The system $\Sigma(\hat{\Gamma})$ is UGS.%
\item[(b)] There exists a path $\sigma$ with properties (i)--(iii) of a $C^0$-path of decay.%
\item[(c)] There exists a path $\sigma$ with properties (i) and (ii) of a $C^0$-path of decay.%
\item[(d)] The set $\Psi(\Gamma)$ is uniformly cofinal.%
\end{enumerate}
\end{proposition}

\begin{proof}
(a) $\Rightarrow$ (b): Let $\varphi \in \KC_{\infty}$ be chosen such that $\|\hat{\Gamma}^n(s)\| \leq \varphi(\|s\|)$. We define $\sigma$ componentwise by%
\begin{equation}\label{eq_def_sigmai}
  \sigma_i(r) := \lim_{n \rightarrow \infty}\hat{\Gamma}^n_i(r\unit) \mbox{\quad for all\ } i \in \IC, \ r \geq 0.%
\end{equation}
The limit exists, because the sequence is increasing and by assumption bounded by $\varphi(\|r\unit\|) = \varphi(r)$. From the construction, it follows that $r\unit \leq \sigma(r) \leq \varphi(r)\unit$ for all $r \geq 0$, showing that $\sigma$ is coercive. Since $\sigma(r)$ is the weak$^*$-limit of the $\Sigma(\hat{\Gamma})$-trajectory starting at $r\unit$, it is a fixed point of $\hat{\Gamma}$, thus a point of decay for $\Gamma$. It follows that $\sigma$ also satisfies the decay property. Finally, $r_1 \leq r_2$ clearly implies $\hat{\Gamma}^n_i(r_1\unit) \leq \hat{\Gamma}^n_i(r_2\unit)$ for all $i \in \IC$ and $n \in \N_0$, and hence $\sigma(r_1) \leq \sigma(r_2)$, showing monotonicity.%

(b) $\Rightarrow$ (c): This is trivial.%

(c) $\Rightarrow$ (d): Assume that $\varphi_1(r)\unit \leq \sigma(r) \leq \varphi_2(r)\unit$ for some $\varphi_1,\varphi_2 \in \KC_{\infty}$. For each $s \in \ell^{\infty}_+(\IC)$, we put $s^* := \sigma(\varphi_1^{-1}(\|s\|))$. It follows that%
\begin{equation*}
  s \leq \|s\|\unit = \varphi_1(\varphi_1^{-1}(\|s\|))\unit \leq \sigma(\varphi_1^{-1}(\|s\|)) = s^*.%
\end{equation*}
Moreover, $\Gamma(s^*) \leq s^*$ by construction and $\|s^*\| \leq \varphi_2 \circ \varphi_1^{-1}(\|s\|)$.%

(d) $\Rightarrow$ (a): Let $\varphi \in \KC_{\infty}$ be such that for every $s \in \ell^{\infty}_+(\IC)$ there is $s^* \in \Psi(\Gamma)$ with $s \leq s^*$ and $\|s^*\| \leq \varphi(\|s\|)$. Then $\hat{\Gamma}^n(s) \leq \hat{\Gamma}^n(s^*) = s^*$ for all $s$ and $n$, implying $\|\hat{\Gamma}^n(s)\| \leq \varphi(\|s\|)$ as required.
\end{proof}

\subsection{Attractivity properties}\label{subsec_stability}

We say that a generalized gain operator $\Gamma$ satisfies the \emph{no-joint-increase (NJI) condition} if%
\begin{equation*}
  \Gamma(s) \not\geq s \mbox{\quad for all\ } s > 0.%
\end{equation*}
This is the classical small-gain condition, which in the classical finite case, with an extra assumption on the MAFs, characterizes the existence of a path of strict decay when required for $\Gamma_{\rho}$ (see Theorem \ref{thm_finite_general_case}). In the general case, the NJI condition is an important ingredient of the following characterization of GATT$^*$.%

\begin{proposition}\label{prop_GATT_star_char}
The following statements are equivalent:%
\begin{enumerate}
\item[(a)] The system $\Sigma(\Gamma)$ is GATT$^{\,*}$.%
\item[(b)] The trajectories of $\Sigma(\Gamma)$ are norm-bounded and $\Gamma$ satisfies the NJI condition.%
\end{enumerate}
\end{proposition}

\begin{proof}
(a) $\Rightarrow$ (b): A trajectory which weakly$^*$ converges is norm-bounded by Lemma \ref{lem_weakstar1}. Hence, every trajectory of $\Sigma(\Gamma)$ is norm-bounded. The NJI condition holds, because otherwise there would exist $s > 0$ with $\Gamma(s) \geq s$ implying $\Gamma^n(s) \geq s$ for all $n$, and hence $\Gamma^n(s)$ could not converge to zero componentwise.%

(b) $\Rightarrow$ (a): For $s \in \ell^{\infty}_+(\IC)$, let $\omega_*(s)$ denote the weak$^*$ $\omega$-limit set of $s$, i.e.~the set of points in $\ell^{\infty}_+(\IC)$ which are weak$^*$-limits of subsequences of $(\Gamma^n(s))_{n\in\N_0}$. By assumption, there exists $b \in \ell^{\infty}_+(\IC)$ such that $\Gamma^n(s)$ is contained in the closed order interval $[0,b]$ for all $n$. Since closed order intervals are sequentially compact in the weak$^*$-topology by Lemma \ref{lem_weakstar3}, $\omega_*(s)$ is nonempty. Now, fix $s \in \ell^{\infty}_+(\IC)$. Then there exists $u \in \ell^{\infty}_+(\IC)$ with $\omega_*(s) \leq u$. From monotonicity and $\Gamma(\omega_*(s)) = \omega_*(s)$,\footnote{This is a general property of $\omega$-limit sets, which also holds in this context, since $\Gamma$ is sequentially continuous in the weak$^*$-topology.} it follows that $\omega_*(s) \leq \Gamma^n(u)$ for all $n$, implying $\omega_*(s) \leq \omega_*(u)$. By the same arguments, there is $v \in \ell^{\infty}_+(\IC)$ with $\omega_*(u) \leq \omega_*(v)$. The set%
\begin{equation*}
  F := \left\{ \tilde{s} \in \ell^{\infty}_+(\IC)\ :\ \omega_*(s) \leq \tilde{s} \leq \omega_*(v) \right\}%
\end{equation*}
is the intersection of closed order intervals. Hence, it is weakly$^*$ sequentially compact and convex. Moreover, it contains $\omega_*(u)$, and hence is nonempty. Since $\Gamma$ is weakly$^*$ sequentially continuous, the set $\Gamma(F)$ is weakly$^*$ sequentially compact. Monotonicity of $\Gamma$ and invariance of $\omega$-limit sets implies $\Gamma(F) \subset F$. It follows from Schauder's fixed point theorem that there is a fixed point in $F$. From the NJI condition, it follows that the origin is the only fixed point. Hence, $0 \in F$, implying $\omega_*(s) = \{0\}$, which means that $\Gamma^n(s) \stackrel{\star}{\rightharpoonup} 0$.
\end{proof}

\begin{remark}
The proof of the preceding proposition is essentially taken from the much more general result \cite[Lem.~1]{Dan} about the existence of fixed points for monotone operators with bounded orbits.
\end{remark}

We now aim at a characterization of GATT. The following lemma will be useful.%

\begin{lemma}\label{lem_aux_stability}
Given a generalized gain operator $\Gamma$, for all $n \in \N$ and $r,\ep > 0$ there exists $\delta > 0$ such that for every $s \geq 0$ with $\|s\| \leq r$ and every $i \in \IC$ the following holds: $\Gamma_j(s) \geq s_j - \delta$ for all $j \in \NC^-_i(n-1)$ implies $\Gamma^n_i(s) \geq s_i - \ep$.
\end{lemma}

\begin{proof}
We put $\delta_n := \ep$ and choose $\ep_i,\delta_i > 0$, $1 \leq i \leq n - 1$, according to the following algorithm: We start with $\ep_{n-1} := \delta_n/2$ and continue recursively by%
\begin{equation*}
  \|s^1 - s^2\| \leq \delta_l \quad \Rightarrow \quad \|\Gamma(s^1) - \Gamma(s^2)\| \leq \ep_l, \quad \ep_{l-1} := \frac{\delta_l}{2}%
\end{equation*}
whenever $\max\{\|s^1\|,\|s^2\|\} \leq r$, which is possible by the uniform continuity of $\Gamma$ on bounded sets. We do this for $l = n - 1$ down to $l = 2$. Eventually, we choose $\delta_1 =: \delta$ satisfying the two requirements $\delta_1 \leq \min\left\{ \ep_i : 1 \leq i \leq n - 1 \right\}$ and%
\begin{equation*}
  \|s^1 - s^2\| \leq \delta_1 \quad \Rightarrow \quad \|\Gamma(s^1) - \Gamma(s^2)\| \leq \ep_1.%
\end{equation*}
Now, let $s \geq 0$ with $\|s\| \leq r$ and $i \in \IC$ such that $\Gamma_j(s) \geq s_j - \delta$ for all $j \in \NC^-_i(n-1)$. Inductively, we get estimates for $\Gamma^l_i(s)$, $1 \leq l \leq n$. First,%
\begin{equation*}
  \Gamma_j(s) \geq s_j - \delta = s_j - \delta_1 \mbox{\quad for all\ } j \in \NC^-_i(n-1)%
\end{equation*}
by assumption. Using, as the induction hypothesis, that $\Gamma^l_j(s) \geq s_j - \delta_l$ for all $j \in \NC^-_i(n - l)$ with $l \in \{1,\ldots,n-1\}$, for each $j \in \NC^-_i(n - l - 1)$ we obtain%
\begin{align*}
  \Gamma^{l + 1}_j(s) &= \Gamma_j(\Gamma^l(s)) = \Gamma_j( [\Gamma^l_k(s)]_{k \in \IC_j} ) \geq \Gamma_j( [\max\{0,s_k - \delta_l\}]_{k \in \IC_j} ) \\
	&\geq \Gamma_j(s) - \ep_l \geq s_j - (\delta_1 + \ep_l) \geq s_j - 2\ep_l = s_j - \delta_{l + 1}.%
\end{align*}
For $l = n - 1$, this yields the desired estimate $\Gamma^n_i(s) \geq s_i - \delta_n = s_i - \ep$.
\end{proof}

We now introduce a uniform version of the NJI condition.%

\begin{definition}
We say that a generalized gain operator $\Gamma$ satisfies the \emph{topologically uniform no-joint-increase (TU-NJI) condition} if the following holds: For all $r,\ep > 0$ there are $n \in \N$ and $\delta > 0$ such that for all $s \in \ell^{\infty}_+(\IC)$ and $i \in \IC$ the following implication holds:%
\begin{equation*}
  s_i \geq \ep  \ \wedge \ \|s\| \leq r  \quad \Rightarrow \quad \exists j \in \NC^-_i(n):\ s_j \geq \delta \ \wedge \ \Gamma_j(s) < s_j.%
\end{equation*}
\end{definition}

The attribute ``topological'' in the name for this property refers to the graph topology of $\GC$. The next proposition shows that in the finite case the TU-NJI condition reduces to the NJI condition.%

\begin{proposition}\label{prop_finite_case_uniform_NJI}
If $|\IC| < \infty$, then the following statements are equivalent:%
\begin{enumerate}
\item[(a)] $\Gamma$ satisfies the NJI condition.%
\item[(b)] $\Gamma$ satisfies the TU-NJI condition.%
\end{enumerate}
\end{proposition}

\begin{proof}
The implication (b) $\Rightarrow$ (a) is easy to see and also holds in the infinite case. To prove the converse implication, we assume that (a) holds, but (b) does not. Then there are $r,\ep > 0$ such that for all $n$ and $\delta$ there exist $s$ and $i$ with $s_i \geq \ep$, $\|s\| \leq r$ and $\Gamma_j(s) \geq s_j$ or $s_j < \delta$ for all $j \in \NC^-_i(n)$. This allows to construct sequences $(s^n)_{n \in \N}$ in $\ell^{\infty}_+(\IC)$ and $(i_n)_{n \in \N}$ in $\IC$ with $s^n_{i_n} \geq \ep$, $\|s^n\| \leq r$, $\Gamma_j(s^n) \geq s^n_j$ or $s^n_j < 1/n$ for all $j \in \NC^-_{i_n}(n)$. Due to the finiteness of $\IC$, we have $\NC^-_{i_n}(n) = \NC^-_{i_n}$ for all sufficiently large $n$ and $i_n \equiv i$ for infinitely many $n$ and a fixed index $i$. Hence, passing to a subsequence, we can assume that $s^n_i \geq \ep$ and $\Gamma_j(s^n) \geq s^n_j$ or $s^n_j < 1/n$ for all $j \in \NC^-_i$. Since $\IC$ is finite, $\ell^{\infty}_+(\IC)$ is finite-dimensional, and hence all $s^n$ are contained in a compact set. Thus, passing to another subsequence if necessary, we can assume that $s^n$ converges to some $s$ with $\|s\| \leq r$ and $s_i \geq \ep$. Since $\Gamma_j(s^n) \geq s^n_j - 1/n$ for all $n$, this $s$ must satisfy $\Gamma_j(s) \geq s_j$ for all $j \in \NC^-_i$. To see that this contradicts (a), we define a new vector $\tilde{s} \in \ell^{\infty}_+(\IC)$ by%
\begin{equation*}
  \tilde{s}_j := \left\{\begin{array}{cc}
                           s_j & \mbox{if } j \in \NC^-_i, \\
													   0 & \mbox{otherwise}.
													\end{array}\right.%
\end{equation*}
Since $s_i > 0$ and $i \in \NC^-_i$, we have $\tilde{s} > 0$. Then (a) implies $\Gamma_j(\tilde{s}) < \tilde{s}_j$ for some index $j$. By construction, we must have $j \in \NC^-_i$, since otherwise $\Gamma_j(\tilde{s})$ would be negative. Since $\Gamma_j(s)$ only depends on those components $s_k$ with $k \in \IC_j \subset \NC^-_i$, it follows that $\Gamma_j(s) = \Gamma_j(\tilde{s}) < \tilde{s}_j = s_j$. The proof is complete.
\end{proof}

The following two propositions relate the TU-NJI condition to GATT.%

\begin{proposition}\label{prop_GATT_implies_uniform_NJI}
If $\Sigma(\Gamma)$ is GATT, then $\Gamma$ satisfies the TU-NJI condition.
\end{proposition}

\begin{proof}
From the proof of Lemma \ref{lem_ugas}, we see that global attractivity together with monotonicity implies uniform global attractivity, i.e.%
\begin{equation*}
  \forall r > 0\ \forall \ep > 0\ \exists n \in \N:\ k \geq n \mbox{ and } \|s\| \leq r \quad \Rightarrow \quad \|\Gamma^k(s)\| \leq \ep.%
\end{equation*}
Now, assume to the contrary that there exist $\ep,r > 0$ such that for every $n \in \N$ and $\delta > 0$ there are $s \in \ell^{\infty}_+(\IC)$ and $i \in \IC$ with $s_i \geq \ep$ and $\|s\| \leq r$, but $\Gamma_j(s) \geq s_j$ or $s_j < \delta$ for all $j \in \NC^-_i(n)$. We then choose $n$ large enough such that $\|s\| \leq r$ implies $\|\Gamma^n(s)\| < \ep/2$, and subsequently, we choose $\delta = \delta(n,r,\ep/2) > 0$ according to Lemma \ref{lem_aux_stability}. By our assumption, there exist $s \geq 0$ and $i \in \IC$ with%
\begin{equation*}
  \ep \leq s_i \leq \Gamma^n_i(s) + \frac{\ep}{2} \leq \|\Gamma^n(s)\| + \frac{\ep}{2} < \frac{\ep}{2} + \frac{\ep}{2} = \ep,%
\end{equation*}
which is a contradiction. To apply the lemma, we have used that both $\Gamma_j(s) \geq s_j$ and $s_j < \delta$ imply $\Gamma_j(s) \geq s_j - \delta$.
\end{proof}

\begin{proposition}\label{prop_uniform_NJI_implies_GATT}
Let the following assumptions hold for a locally finite generalized gain operator $\Gamma$ and some $\rho \in \KC_{\infty}$:%
\begin{enumerate}
\item[(i)] The set $\Psi(\Gamma_{\rho})$ is cofinal.%
\item[(ii)] $\Gamma_{\rho}$ satisfies the TU-NJI condition.%
\item[(iii)] There exists $M > 0$ with $|\IC_i| \leq M$ for all $i \in \IC$.%
\end{enumerate}
Then there exists $\rho' \in \KC_{\infty}$ such that $\Sigma(\Gamma_{\rho'})$ is GATT. 
\end{proposition}

\begin{proof}
For notational simplicity, we write the proof for $\rho' = 0$, i.e.~we show that $\Sigma(\Gamma)$ is GATT. We can do this without loss of generality, since we can decompose $\id + \rho = (\id + \rho_1) \circ (\id + \rho_2)$ according to Lemma \ref{lem_kinfty} and use $\rho' = \rho_1$, which replaces $\Gamma$ with $\Gamma_{\rho_1}$ throughout the proof. Observe that $\Psi(\Gamma_{\rho}) \subset \Psi(\Gamma)$, hence assumption (i) implies that $\Psi(\Gamma)$ is cofinal. It then suffices to prove that every trajectory starting in $\Psi(\Gamma)$ converges to zero. To this end, fix $s \in \Psi(\Gamma)$ and $\ep > 0$ with $\ep\unit \leq s$, and consider the decreasing trajectory $s^n := \Gamma^n(s)$, $n \in \N_0$. Note that $s^n \stackrel{\star}{\rightharpoonup} 0$, because assumption (ii) implies that the origin is the only fixed point of $\Gamma$.%

By assumption (ii), there are $m = m(\ep,\|s\|) \in \N$ and $\delta = \delta(\ep,\|s\|) > 0$ such that for every $i \in \IC$ we have $\Gamma_j(s^n) < (\id + \rho)^{-1}(s^n_j)$ for some $j = j(n) \in \NC^-_i(m)$ with $s_j^n \geq \delta$ as long as $s^n_i \geq \ep$. We define%
\begin{equation*}
  K_i := \max\{ k \in \N : s^k_i \geq \ep \} \mbox{\quad for all\ } i \in \IC.%
\end{equation*}
Then, for each $i \in \IC$ there must exist some $j_i \in \NC^-_i(m)$ with%
\begin{equation*}
  |\{ k \in [0;K_i-1] : s^k_{j_i} \geq \delta \wedge \Gamma_{j_i}(s^k) < (\id + \rho)^{-1}(s^k_{j_i}) \}| \geq \frac{K_i}{|\NC^-_i(m)|}.%
\end{equation*}
If this was not the case, we would end up with the contradiction%
\begin{align*}
  K_i &= |\{ k \in [0;K_i-1] : \exists j \in \NC^-_i(m),\	s^k_j \geq \delta \wedge \Gamma_j(s^k) < (\id + \rho)^{-1}(s^k_j) \}| \\
	    &\leq \sum_{j \in \NC^-_i(m)} |\{ k \in [0;K_i-1] : s^k_j \geq \delta \wedge \Gamma_j(s^k) < (\id + \rho)^{-1}(s^k_j) \}| < K_i.%
\end{align*}
From assumption (iii), it follows that $|\NC^-_i(m)| \leq B_m$ with $B_m := \sum_{l=0}^m M^l$. Now, let $k_i \in [0;K_i - 1]$ be the largest time at which $s^k_{j_i} \geq \delta$ and $\Gamma_{j_i}(s^k) < (\id + \rho)^{-1}(s^k_{j_i})$. Then we obtain%
\begin{equation*}
  \delta \leq s^{k_i}_{j_i} \leq (\id + \rho)^{-(\lceil K_i/|\NC^-_i(m)| \rceil - 1)}(s_{j_i}) \leq (\id + \rho)^{-(\lceil K_i/B_m \rceil - 1)}(\|s\|).%
\end{equation*}
This implies that there is a bound on $K_i$, which is uniform in $i$, since the right-hand side converges to zero as $K_i \rightarrow \infty$. But if $K_i \leq K$ for all $i \in \IC$, then $\|\Gamma^{K+1}(s)\| < \ep$. The number $\ep > 0$ can be chosen arbitrarily small. Hence, we have shown that $\Gamma^n(s)$ converges to zero, which completes the proof.
\end{proof}

\section{Existence of paths of strict decay}\label{sec_path_of_decay}

In this section, we derive necessary conditions and sufficient conditions for the existence of a path of strict decay. Recall that by Proposition \ref{prop_path_of_decay_2} it suffices to consider $C^0$-paths.

\subsection{The \texorpdfstring{\(\oplus\)}{max}-MBI property}\label{subsec_max_mbi}

We first introduce a property that will turn out to be crucial for characterizing the existence of a path of strict decay. This property first appeared in \cite{Ruf}, but was only used for max-type gain operators there.

\begin{definition}
A generalized gain operator $\Gamma$ satisfies the \emph{$\oplus$-MBI property} if there exists $\varphi \in \KC_{\infty}$ such that for all $s,b \in \ell^{\infty}_+(\IC)$ the following implication holds:%
\begin{equation*}
  s \leq \bsign{\Gamma}_b(s) \quad \Rightarrow \quad \|s\| \leq \varphi(\|b\|).%
\end{equation*}
\end{definition}

\begin{remark}
The letters \emph{MBI} stand for \emph{monotone bounded invertibility}. As an explanation for this name, think of $\Gamma$ being a linear operator, $\oplus$ being replaced with addition, and $\varphi$ a linear function. Then the required property could be rewritten as%
\begin{equation*}
  (\id - \Gamma)s \leq b \quad \Rightarrow \quad \|s\| \leq K\|b\|,%
\end{equation*}
with a constant $K > 0$. If we assume that $\id - \Gamma$ is an invertible bounded operator with positive inverse, this property is satisfied with $K = \|(\id - \Gamma)^{-1}\|$. Hence, the $\oplus$-MBI property generalizes the bounded invertibility and positivity of the operator $\id - \Gamma$. By \cite[Thm.~3.3]{GMi}, this property is equivalent to the spectral radius of $\Gamma$ being strictly less than one. Also recall that the \emph{bounded inverse theorem} from functional analysis guarantees that an invertible bounded operator between Banach spaces always has a bounded inverse.
\end{remark}

\begin{remark}
It is easy to see that for verifying the $\oplus$-MBI property one only needs to consider vectors of the form $b = r\unit$ with $r \geq 0$.
\end{remark}

Although we will not use the following proposition, we present it to demonstrate that the $\oplus$-MBI property is a natural property. For this result, observe that iterates of generalized gain operators are generalized gain operators as well.

\begin{proposition}
If a generalized gain operator $\Gamma$ satisfies the $\oplus$-MBI property, then so does any iterate of $\Gamma$, with the same $\KC_{\infty}$-function $\varphi$.
\end{proposition}

\begin{proof}
Assume that $\Gamma$ satisfies the $\oplus$-MBI property and let $s \leq b \oplus \Gamma^m(s)$ for some $s,b \geq 0$ and $m \in \N$. We put%
\begin{equation}\label{eq_def_tildes}
  \tilde{s} := \bigoplus_{l=0}^{m-1} \Gamma^l(s) = s \oplus \Gamma(s) \oplus \cdots \oplus \Gamma^{m-1}(s).%
\end{equation}
Due to the monotonicity of $\Gamma$, we have%
\begin{equation*}
  \Gamma(\tilde{s}) \geq \bigoplus_{l=0}^{m-1} \Gamma^{l + 1}(s),%
\end{equation*}
which leads to%
\begin{align*}
  b \oplus \Gamma(\tilde{s}) &\geq b \oplus \bigoplus_{l=0}^{m-1} \Gamma^{l + 1}(s) = b \oplus \Gamma(s) \oplus \Gamma^2(s) \oplus \cdots \oplus \Gamma^m(s) \\
																&= \left[ b \oplus \Gamma^m(s) \right] \oplus \left[ \Gamma(s) \oplus \Gamma^2(s) \oplus \cdots \oplus \Gamma^{m-1}(s) \right] \\
																&\geq s \oplus \Gamma(s) \oplus \cdots \oplus \Gamma^{m-1}(s) = \tilde{s}.%
\end{align*}
Hence, $\|\tilde{s}\| \leq \varphi(\|b\|)$, implying $\|s\| \leq \|\tilde{s}\| \leq \varphi(\|b\|)$.
\end{proof}

The next proposition provides several characterizations.%

\begin{proposition}\label{prop_maxmbi_char}
For a generalized gain operator $\Gamma$, let $A := \bigcup_{r \geq 0}\Fix(\bsign{\Gamma}_{r\unit}) \subset \Psi(\Gamma)$. Then the following statements are equivalent:%
\begin{enumerate}
\item[(a)] $\Gamma$ satisfies the $\oplus$-MBI property.%
\item[(b)] The system $\Sigma(\hat{\Gamma})$ is UGS, the system $\Sigma(\Gamma)$ is GATT$^{\,*}$ and $A$ is coercive.%
\item[(c)] The system $\Sigma(\hat{\Gamma})$ is UGS, $\Gamma$ satisfies the NJI condition and $A$ is coercive.%
\item[(d)] The set $\Psi(\Gamma)$ is cofinal, $\Gamma$ satisfies the NJI condition and $A$ is coercive.%
\end{enumerate}
\end{proposition}

\begin{proof}
(a) $\Rightarrow$ (b): We first show that $\Sigma(\hat{\Gamma})$ is UGS. By Proposition \ref{prop_hatgamma_properties} (f), we have $\hat{\Gamma}^{n+1}(s) \equiv s \oplus \Gamma(\hat{\Gamma}^n(s))$, implying $\hat{\Gamma}^n(s) \leq s \oplus \Gamma(\hat{\Gamma}^n(s))$. The $\oplus$-MBI property thus immediately implies $\|\hat{\Gamma}^n(s)\| \leq \varphi(\|s\|)$ for an appropriate $\varphi \in \KC_{\infty}$, showing that $\Sigma(\hat{\Gamma})$ is UGS. Now, we prove that $\Sigma(\Gamma)$ is GATT$^*$ via Proposition \ref{prop_GATT_star_char}. Putting $b = 0$, we see that the $\oplus$-MBI property implies the NJI condition. Moreover, since $\Sigma(\hat{\Gamma})$ is UGS, also $\Sigma(\Gamma)$ is UGS, thus the trajectories of $\Sigma(\Gamma)$ are norm-bounded. Finally, we show that $A$ is coercive. Let $s = r\unit \oplus \Gamma(s)$. By assumption, $\|s\| \leq \varphi(r)$ for some $\varphi \in \KC_{\infty}$, hence $s \geq r\unit \geq \varphi^{-1}(\|s\|)\unit$.%

(b) $\Rightarrow$ (c): This follows from Proposition \ref{prop_GATT_star_char}.%

(c) $\Rightarrow$ (d): Obviously, $\Psi(\Gamma)$ is cofinal if $\Sigma(\hat{\Gamma})$ is UGS.%

(d) $\Rightarrow$ (a): Let $s \leq \bsign{\Gamma}_{r\unit}(s)$ for some $s \in \ell^{\infty}_+(\IC)$ and $r \geq 0$.  Putting $s^n := \bsign{\Gamma}_{r\unit}^n(s)$, this implies $s^n \leq s^{n+1}$ for all $n$. By assumption, there exists $\hat{s} \in \Psi(\Gamma)$ with $\hat{s} \geq r\unit \oplus s$. This implies $s^n \leq \bsign{\Gamma}_{r\unit}^n(\hat{s})$ for all $n$ and $\bsign{\Gamma}_{r\unit}(\hat{s}) = r\unit \oplus \Gamma(\hat{s}) \leq r\unit \oplus \hat{s} = \hat{s}$ so that $s^n \leq \hat{s}$ for all $n$. Hence, $s^n$ weakly$^*$ converges to a fixed point $s^*$ of $\bsign{\Gamma}_{r\unit}$. By assumption, $s^* \geq \underline{\varphi}(\|s^*\|)\unit$, implying $\|s^*\| \leq \underline{\varphi}^{-1}(\inf_{i\in\IC}s^*_i)$. Moreover, there exists $i \in \IC$ with $\Gamma_i(s^*) < s^*_i$, implying $s^*_i = r$. Altogether, $\|s\| \leq \|s^*\| \leq \underline{\varphi}^{-1}(s^*_i) = \underline{\varphi}^{-1}(r)$.
\end{proof}

\begin{corollary}\label{cor_unique_fp}
Under the additional assumption that each of the operators $\bsign{\Gamma}_b$ has a unique fixed point, the statements (a)--(c) in Proposition \ref{prop_maxmbi_char} are still equivalent after removing the requirement that $A$ is coercive.
\end{corollary}

\begin{proof}
It suffices to prove that $\Sigma(\hat{\Gamma})$ being UGS together with the additional assumption implies that $A$ is coercive. To this end, let $s^*$ be the unique fixed point of $\bsign{\Gamma}_{r\unit}$ for a fixed $r \geq 0$. The trajectory $s^n := \bsign{\Gamma}_{r\unit}^n(r\unit)$ is increasing and since $r\unit \leq s^*$, we have $s^n \stackrel{\star}{\rightharpoonup} s^*$. Recall that $\hat{\Gamma}^n(r\unit) \equiv s^n$, and thus $\|s^n\| \leq \varphi(r)$ for a $\KC_{\infty}$-function $\varphi$, independent of $r$. This implies $r\unit \leq s^* \leq \varphi(r)\unit$, and thus $s^* \geq \varphi^{-1}(\|s^*\|)\unit$, showing the coercivity of $A$.
\end{proof}

Next, we show that the $\oplus$-MBI property implies that each of the operators $\bsign{\Gamma}_b$ has a minimal and a maximal fixed point.%

\begin{proposition}\label{prop_min_max_fixed_points}
Let $\Gamma$ be a generalized gain operator satisfying the $\oplus$-MBI property. Then, for each $b \in \ell^{\infty}_+(\IC)$, there exist fixed points $s_*(b)$ and $s^*(b)$ of $\bsign{\Gamma}_b$ such that the following statements hold:%
\begin{enumerate}
\item[(a)] Every $s \in \Fix(\bsign{\Gamma}_b)$ satisfies $s_*(b) \leq s \leq s^*(b)$.%
\item[(b)] If $b^1 \leq b^2$, then $s_*(b_1) \leq s_*(b_2)$ and $s^*(b_1) \leq s^*(b_2)$.%
\item[(c)] For every $s \leq s_*(b)$, we have $\bsign{\Gamma}_b^n(s) \stackrel{\star}{\rightharpoonup} s_*(b)$.%
\item[(d)] For every $s \geq s^*(b)$, we have $\bsign{\Gamma}_b^n(s) \stackrel{\star}{\rightharpoonup} s^*(b)$.
\end{enumerate}
\end{proposition}

\begin{proof}
We define $s_*(b)$ as the weak$^*$ limit of the trajectory $s^n := \bsign{\Gamma}_b^n(b)$, $n \in \N_0$. Since $\bsign{\Gamma}_b(b) \geq b$ and $\bsign{\Gamma}_b$ is monotone, this trajectory is increasing. Since $s^n \leq s^{n+1} = b \oplus \Gamma(s^n)$ implies $\|s^n\| \leq \varphi(\|b\|)$ for all $n \in \N_0$ by the $\oplus$-MBI property, the trajectory is also bounded. Hence, the weak$^*$ limit exists and must be a fixed point of $\bsign{\Gamma}_b$. If $s$ is an arbitrary fixed point of $\bsign{\Gamma}_b$, we must have $s \geq b$, and thus $s = \bsign{\Gamma}_b^n(s) \geq s^n$ for all $n$, implying $s \geq s_*(b)$. This shows that $s_*(b)$ is the minimal fixed point of $\bsign{\Gamma}_b$.%

We now prove the existence of the maximal fixed point $s^*(b)$. To this end, we put $r := \varphi(\|b\|)$ and consider the point $s^0 := s_*(r\unit)$ and its trajectory $s^n := \bsign{\Gamma}_b^n(s^0)$, $n \in \N_0$. We have%
\begin{equation*}
  s^1 = b \oplus \Gamma(s_*(r\unit)) \leq r\unit \oplus \Gamma(s_*(r\unit)) = s_*(r\unit) = s^0.%
\end{equation*}
By monotonicity of $\bsign{\Gamma}_b$, this implies that the trajectory $(s^n)_{n\in\N_0}$ is decreasing, and thus $s^n \stackrel{\star}{\rightharpoonup} s^*$ for some $s^* \in \Fix(\bsign{\Gamma}_b)$. If $s$ is an arbitrary fixed point of $\bsign{\Gamma}_b$, the $\oplus$-MBI property implies $s \leq \varphi(\|b\|)\unit \leq s_*(\varphi(\|b\|)\unit) = s^0$, so that $s = \bsign{\Gamma}_b^n(s) \leq \bsign{\Gamma}_b^n(s^0) = s^n$ for all $n \in \N_0$. Hence, $s \leq s^*$, showing that $s^* = s^*(b)$ is the desired maximal fixed point.%

It remains to prove (b)--(d). Statement (b) is an easy consequence of monotonicity. Every trajectory starting in $s \leq s_*(b)$ is bounded below by the increasing trajectory starting in the origin, which implies (c). Conversely, every trajectory starting in $s \geq s^*(b)$ is bounded by a decreasing trajectory as constructed above, implying (d).
\end{proof}

\begin{figure}[ht!]\label{fig:gamma_b_dynamics}
 \centering
  \includegraphics[width=0.65\textwidth]{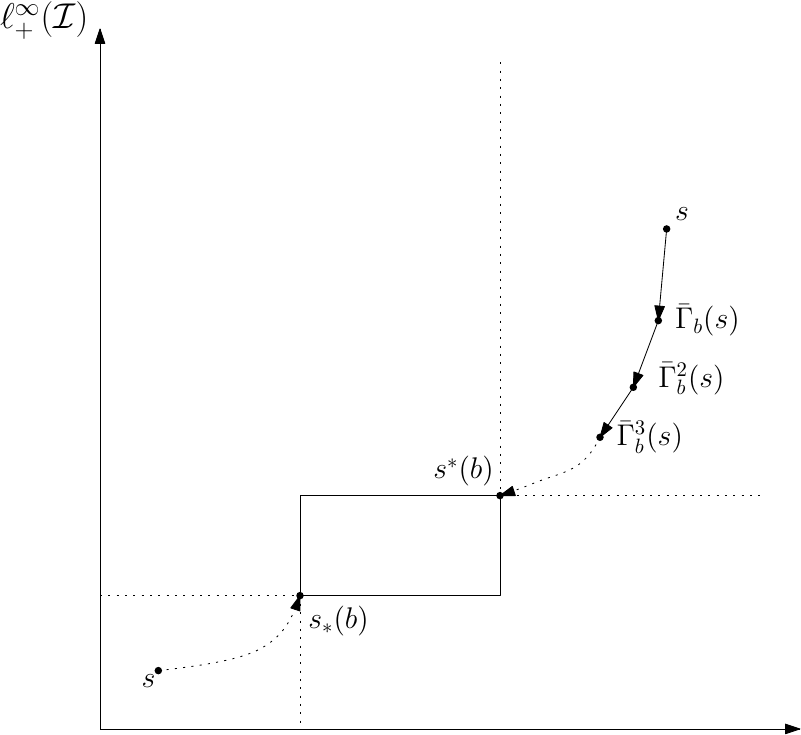}
 \caption{The dynamics of $\bsign{\Gamma}_b$ in the weak$^*$-topology under the $\oplus$-MBI property. The order interval $[s_*(b),s^*(b)]$ is an invariant set which is a weak$^*$ attractor. (Observe that any $s$ not satisfying $s \leq s_*(b)$ or $s \geq s^*(b)$ does satisfy $s^1 \leq s \leq s^2$ for some $s^1 \leq s_*(b)$ and $s^*(b) \leq s^2$.)}
\end{figure}

Figure \ref{fig:gamma_b_dynamics} illustrates the dynamics of the projected gain operator $\bsign{\Gamma}_b$ in the weak$^*$-topology.%

Finally, we derive a characterization of the $\oplus$-MBI property in terms of the componentwise decay in one time step, similar to the TU-NJI condition.%

\begin{proposition}\label{prop_maxmbi_nji_char}
For a generalized gain operator $\Gamma$, the following statements are equivalent:%
\begin{enumerate}
\item[(a)] $\Gamma$ satisfies the $\oplus$-MBI property.%
\item[(b)] There exists $\varphi \in \KC_{\infty}$ such that for all $s \in \ell^{\infty}_+(\IC)$ and $i \in \IC$ with $s_i > 0$ there exists $j \in \NC^-_i$ with $s_j \geq \varphi(s_i)$ and $\Gamma_j(s) < s_j$.
\end{enumerate}
\end{proposition}

\begin{proof}
(a) $\Rightarrow$ (b): Given $s \in \ell^{\infty}_+(\IC)$, we define the vector $b$ by%
\begin{equation*}
  b_i := \left\{\begin{array}{cc} 
	                  0 & \mbox{if } \Gamma_i(s) \geq s_i, \\
										s_i & \mbox{if } \Gamma_i(s) < s_i.
							  \end{array}\right.%
\end{equation*}
By definition, it then holds that $s \leq b \oplus \Gamma(s)$. Thus, by assumption, we have%
\begin{equation*}
  \sup_{i\in\IC}s_i = \|s\| \leq \varphi(\|b\|) = \sup\left\{ \varphi(s_j)\ :\ j \in \IC,\ \Gamma_j(s) < s_j \right\}.%
\end{equation*}
Let us fix an index $i \in \IC$ with $s_i > 0$. Since (a) implies the NJI condition, the set $\{ j \in \IC\ :\ \Gamma_j(s) < s_j \}$ is nonempty. If we put $\tilde{\varphi} := 2\varphi$, we then have%
\begin{equation*}
  s_i < \sup\left\{ \tilde{\varphi}(s_j)\ :\ \Gamma_j(s) < s_j \right\}.%
\end{equation*}
This implies the existence of $j \in \IC$ with $s_j \geq \tilde{\varphi}^{-1}(s_i)$ and $\Gamma_j(s) < s_j$. It remains to prove that $j$ can be chosen as an element of $\NC^-_i$. To this end, define the vector $\tilde{s}$ by%
\begin{equation*}
  \tilde{s}_j := \left\{ \begin{array}{cc}
	                                s_j & \mbox{if } j \in \NC^-_i, \\
																	  0 & \mbox{otherwise}.
													\end{array}\right.%
\end{equation*}
Since $\tilde{s}_i = s_i > 0$, there is $j \in \IC$ with $\Gamma_j(\tilde{s}) < \tilde{s}_j$ and $\tilde{s}_j \geq \tilde{\varphi}^{-1}(\tilde{s}_i)$ (by what we have shown already). By definition of $\tilde{s}$, it follows that $j \in \NC^-_i$ and thus $\tilde{s}_j = s_j$. Moreover, $\Gamma_j(\tilde{s}) = \Gamma_j(s)$, since $\Gamma_j(\tilde{s})$ only depends on the components $\tilde{s}_k$ with $k \in \IC_j \subset \NC^-_i$.%

(b) $\Rightarrow$ (a): Let $s \leq b \oplus \Gamma(s)$. By assumption, there exists $\varphi$ such that for every $i \in \IC$ with $s_i > 0$ there is $j \in \NC^-_i$ with $s_j \geq \varphi(s_i)$ and $\Gamma_j(s) < s_j$. The inequality $\Gamma_j(s) < s_j$ in combination with $s_j \leq \max\{b_j,\Gamma_j(s)\}$ implies $s_j \leq b_j$, and hence $s_i \leq \varphi^{-1}(b_j) \leq \varphi^{-1}(\|b\|)$. Altogether, $\|s\| \leq \varphi^{-1}(\|b\|)$.
\end{proof}

The proposition leads to the following definition.%

\begin{definition}
We say that a generalized gain operator $\Gamma$ satisfies the \emph{spatially uniform no-joint-increase (SU-NJI) condition} if the following holds: There exists $\varphi \in \KC_{\infty}$ such that for every $s \in \ell^{\infty}_+(\IC)$ and $i \in \IC$ the following implication holds:%
\begin{equation*}
  s_i > 0 \quad \Rightarrow \quad \exists j \in \NC^-_i:\ s_j \geq \varphi(s_i) \ \wedge \ \Gamma_j(s) < s_j.%
\end{equation*}
\end{definition}

\begin{remark}
The TU-NJI and the SU-NJI condition have a similar structure. Both are uniform versions of the NJI condition. The uniformity in the TU-NJI condition concerns the distance of the vertices $i$ and $j$ in the graph $\GC$, while the uniformity in the SU-NJI condition concerns the relation between $s_i$ and $s_j$. Indeed, the plain NJI condition can be written as%
\begin{equation*}
  \forall s \in \ell^{\infty}_+(\IC),\ \forall i \in \IC:\quad s_i > 0 \quad \Rightarrow \quad \exists j \in \NC^-_i:\ \Gamma_j(s) < s_j.%
\end{equation*}
Hence, the NJI condition and the SU-NJI condition only differ in the estimate $s_j \geq \varphi(s_i)$. In the finite case, it follows from Proposition \ref{prop_finite_case_uniform_NJI} that the SU-NJI condition implies the TU-NJI condition. In the infinite case, however, it is quite clear that the SU-NJI condition does not imply the TU-NJI condition, since the uniformity in the distance between the vertices $i$ and $j$ is completely missing in the first.
\end{remark}

\begin{remark}\label{rem_uniform_sgc}
The SU-NJI condition has a strong similarity to the \emph{uniform small-gain condition} from \cite{KMG}: There exists $\varphi \in \KC_{\infty}$ with%
\begin{equation*}
  \|s\| \leq \sup\left\{ \varphi(s_i - \Gamma_i(s))\ :\ i \in \IC,\ \Gamma_i(s) < s_i \right\} \mbox{\quad for all\ } s \in \ell^{\infty}_+(\IC)%
\end{equation*}
This condition is equivalent to the MBI property, which reads%
\begin{equation*}
  s \leq b + \Gamma(s) \quad \Rightarrow \quad \|s\| \leq \varphi(\|b\|).%
\end{equation*}
In our paper, we do not use these conditions, however.
\end{remark}

\subsection{Necessary conditions}\label{subsec_necessary}

From Proposition \ref{prop_hatsigma_ugs}, we know that UGS of $\Sigma(\hat{\Gamma}_{\rho})$ for some $\rho \in \KC_{\infty}$ is a necessary condition for the existence of a path of strict decay. The following theorem provides further necessary conditions that we obtain from the continuity of the path which we have not used so far.%

\begin{theorem}\label{thm_necessary_conditions}
If a generalized gain operator $\Gamma$ admits a $C^0$-path $\sigma$ of strict decay, the following statements hold for some $\rho \in \KC_{\infty}$:%
\begin{enumerate}
\item[(a)] $\Sigma(\Gamma_{\rho})$ is UGAS.%
\item[(b)] $\Gamma_{\rho}$ satisfies the $\oplus$-MBI property.%
\end{enumerate}
These statements also hold if $\sigma$ only satisfies properties (i), (ii) and (iv) of a $C^0$-path of strict decay, i.e.~it is not necessary that $\sigma$ is increasing.
\end{theorem}

\begin{proof}
We let $\tilde{\rho}$ denote the $\KC_{\infty}$-function such that $\Gamma_{\tilde{\rho}} \circ \sigma \leq \sigma$. We write $\id + \tilde{\rho} = (\id + \rho') \circ (\id + \rho)$ with $\rho,\rho' \in \KC_{\infty}$, using Lemma \ref{lem_kinfty}. It is then easy to see that $\Gamma_{\rho}(\sigma(r)) \leq (\id + \rho')^{-1} \circ \sigma(r)$ for all $r \geq 0$. We further write $(\id + \rho')^{-1} = \id - \eta$ with $\eta \in \KC_{\infty}$, also using Lemma \ref{lem_kinfty}. Then%
\begin{equation}\label{eq_decrease_property}
  \Gamma_{\rho}(\sigma(r)) \leq (\id - \eta)(\sigma(r)) \mbox{\quad for all\ } r \geq 0.%
\end{equation}
By Proposition \ref{prop_path_of_decay} (b) (which does not require monotonicity), we can assume that for some $\varphi \in \KC_{\infty}$ it holds that%
\begin{equation}\label{eq_sigma_bounds}
  r\unit \leq \sigma(r) \leq \varphi(r)\unit \mbox{\quad for all\ } r \geq 0.%
\end{equation}

(a) By Proposition \ref{prop_hatsigma_ugs} (which does not use monotonicity), $\Sigma(\hat{\Gamma}_{\tilde{\rho}})$ is UGS, and hence also $\Sigma(\hat{\Gamma}_{\rho})$ and $\Sigma(\Gamma_{\rho})$.\footnote{In general, if $\Gamma_1$ and $\Gamma_2$ are monotone operators with $\Gamma_1(s) \leq \Gamma_2(s)$ for all $s \geq 0$ and $\Sigma(\Gamma_2)$ is UGS, then so is $\Sigma(\Gamma_1)$.} By Lemma \ref{lem_ugas}, it remains to prove global attractivity of $\Sigma(\Gamma_{\rho})$. To this end, due to coercivity of $\sigma$ it suffices to prove that every trajectory of the form $(\Gamma_{\rho}^n(\sigma(r)))_{n \in \N_0}$ with $r > 0$ converges to the origin. Hence, fix some $r > 0$ and put $s^n := \Gamma_{\rho}^n(\sigma(r))$ for all $n \in \N_0$. We define%
\begin{equation*}
  r_n := \min\{ \tilde{r} \geq 0 : s^n \leq \sigma(\tilde{r}) \}%
\end{equation*}
for every $n \in \N_0$, which exists by continuity and coercivity of $\sigma$. Since $\Gamma_{\rho}(s^0) \leq s^0$ by \eqref{eq_decrease_property}, the sequence $(s^n)_{n\in\N_0}$ is decreasing. Hence, $s^{n+1} \leq s^n \leq \sigma(r_n)$, so that $r_n \in \{\tilde{r} \geq 0 : s^{n+1} \leq \sigma(\tilde{r}) \}$, implying $r_{n+1} \leq r_n$. Consequently, $(r_n)_{n\in\N_0}$ converges to some $r_* \in \R_+$. If $r_* = 0$, then $\|s^n\| \leq \|\sigma(r_n)\| \leq \varphi(r_n) \rightarrow 0$, proving the assertion. Hence, let us assume to the contrary that $r_* > 0$. Let $\ep := \eta(r_*)$. By \eqref{eq_decrease_property} and \eqref{eq_sigma_bounds}, for any $r \geq r_*$, we have%
\begin{equation*}
  \Gamma_{\rho}(\sigma(r)) \leq \sigma(r) - \eta(\sigma(r)) \leq \sigma(r) - \eta(r_*)\unit = \sigma(r) - \ep\unit.%
\end{equation*}
Since $\sigma(r_n)$ converges to $\sigma(r_*)$ by continuity of $\sigma$, we can choose $n_0$ large enough so that $\|\sigma(r_n) - \sigma(r_{n+1})\| \leq \ep/2$ for all $n \geq n_0$, implying%
\begin{equation*}
  s^{n+1} = \Gamma_{\rho}(s^n) \leq \Gamma_{\rho}(\sigma(r_n)) \leq \sigma(r_n) - \ep\unit \leq \sigma(r_{n+1}) - \frac{\ep}{2}\unit.%
\end{equation*}
However, the inequality $s^{n+1} + (\ep/2)\unit \leq \sigma(r_{n+1})$ contradicts the minimality of $r_{n+1}$. Indeed, for any $r < r_{n+1}$ we must have $s^{n+1}_i > \sigma_i(r)$ for some $i \in \IC$. Thus, $\|\sigma(r) - \sigma(r_{n+1})\| \geq \sigma_i(r_{n+1}) - \sigma_i(r) > \ep/2$, which contradicts the continuity of $\sigma$ for $r$ close enough to $r_{n+1}$.%

(b) Let $s \in \ell^{\infty}_+(\IC)$ and $r \geq 0$ be chosen such that $s \leq r\unit \oplus \Gamma_{\rho}(s)$. The case $s = 0$ is trivial, so we assume that $s > 0$. We define%
\begin{equation*}
  r' := \min\{ \tilde{r} \geq 0 : s \leq \sigma(\tilde{r}) \} > 0.%
\end{equation*}
Again, the minimum exists due to continuity. From the monotonicity of $\bsign{\Gamma}_{\rho,r\unit}$ and $s \leq \sigma(r')$, it follows that%
\begin{equation}\label{eq_some_estimates}
  s \leq \bsign{\Gamma}_{\rho,r\unit}(s) \leq \bsign{\Gamma}_{\rho,r\unit}(\sigma(r')) = r\unit \oplus \Gamma_{\rho}(\sigma(r')) \leq r\unit \oplus (\id - \eta)(\sigma(r')).%
\end{equation}
We pick $\ep \in (0,\eta(r'))$. Since $r'$ was chosen minimal, there exists an index $i \in \IC$ with $\sigma_i(r') - s_i < \ep$. Evaluating \eqref{eq_some_estimates} at this index yields%
\begin{equation*}
  \sigma_i(r') - \ep < s_i \leq \max\{r, \sigma_i(r') - \eta(\sigma_i(r')) \} \leq \max\{r, \sigma_i(r') - \ep\}.%
\end{equation*}
To avoid a contradiction, the maximum at the end of the line must be attained at $r$. Since $r' \leq \sigma_i(r')$, this implies $r' - \ep \leq r$. Sending $\ep$ to zero, we obtain $r' \leq r$, and thus $s \leq \sigma(r') \leq \varphi(r')\unit \leq \varphi(r)\unit$, which implies $\|s\| \leq \varphi(r)$. Hence, $\Gamma_{\rho}$ satisfies the $\oplus$-MBI property.%
\end{proof}

\begin{remark}
In general, it is not clear how the two properties in (a) and (b) are related to each other. In some cases, as we will show, one implies the other, but the general case is open. What we know for sure is that the $\oplus$-MBI property does not imply UGAS in the general infinite case. This can be seen by looking at max-type gain operators, where the $\oplus$-MBI property is equivalent to UGS plus GATT$^{\,*}$, which is strictly weaker than UGAS in the general infinite case.
\end{remark}

\subsection{Path construction}\label{subsec_path_construction}

In this subsection, we introduce a method for the construction of a $C^0$-path of strict decay. We assume that $\Gamma$ satisfies the $\oplus$-MBI property, which, in particular, implies that minimal and maximal fixed points $\sigma_*(r) := s_*(r\unit)$ and $\sigma^*(r) := s^*(r\unit)$ of the operator $\bsign{\Gamma}_{r\unit}$ exist for all $r \geq 0$, see Proposition \ref{prop_min_max_fixed_points}. We fix a function $\varphi \in \KC_{\infty}$ such that%
\begin{equation*}
  s \leq \bsign{\Gamma}_b(s) \quad \Rightarrow \quad \|s\| \leq \varphi(\|b\|) \mbox{\quad for all\ } s,b \in \ell^{\infty}_+(\IC).%
\end{equation*}
The construction starts by fixing a strictly increasing double-sided sequence of real numbers $(r_k)_{k\in\Z}$ with%
\begin{equation*}
  \lim_{k \rightarrow -\infty}r_k = 0 \mbox{\quad and \quad} \lim_{k \rightarrow \infty}r_k = \infty.%
\end{equation*}
We also consider the associated sequence of points $s^k := \sigma^*(r_k)$, which has the following properties for all $k \in \Z$:%
\begin{itemize}
\item $s^k \leq s^{k+1}$.
\item $r_k\unit \leq s^k \leq \varphi(r_k)\unit$.
\item $s^k \in \Psi(\Gamma)$, since $\Fix(\bsign{\Gamma}_{r\unit}) \subset \Psi(\bsign{\Gamma}_{r\unit}) \subset \Psi(\Gamma)$.
\end{itemize}
Our goal is to construct a $C^0$-path of decay $\sigma:\R_+ \rightarrow \ell^{\infty}_+(\IC)$ which satisfies $\sigma(r_k) = s^k$ for all $k \in \Z$. Hence, we have to define $\sigma$ on the open intervals $(r_k,r_{k+1})$ and at $r = 0$. Here, we use statements (b) and (d) of Proposition \ref{prop_min_max_fixed_points}, which show that $s^{k,n} := \bsign{\Gamma}_{r_k\unit}^n(s^{k+1}) \stackrel{\star}{\rightharpoonup} s^k$. We define the path on $(r_k,r_{k+1})$ by linear interpolation between the points $s^{k+1} = s^{k,0}$, $s^{k,1}$, $s^{k,2}$, and so on, and put $\sigma(0) := 0$. In this way, we obtain a path $\sigma:\R_+ \rightarrow \ell^{\infty}_+(\IC)$, satisfying the properties listed in the following proposition. Note that the parametrization of the path on the intervals $(r_k,r_{k+1})$ is not unique.%

\begin{proposition}\label{prop_maxmbi_path}
The path $\sigma$ satisfies properties (i)--(iii) of a $C^0$-path of decay. Additionally, it has the following continuity properties:%
\begin{enumerate}
\item[(C1)] It is norm-continuous from below, i.e., if $(\rho_n)_{n \in \N}$ is a sequence in $\R_+$ converging to some $\rho$ from below, then $\|\sigma(\rho_n) - \sigma(\rho)\| \rightarrow 0$.%
\item[(C2)] It is sequentially weakly$^*$ continuous from above, i.e., if $(\rho_n)_{n \in \N}$ is a sequence in $\R_+$ converging to some $\rho$ from above, then $\sigma(\rho_n) \stackrel{\star}{\rightharpoonup} \sigma(\rho)$.%
\end{enumerate}
\end{proposition}

\begin{proof}
Property (i) (decay): The points $s^k = \sigma(r_k)$ are points of decay of $\bsign{\Gamma}_{r_k\unit}$, hence also points of decay of $\Gamma$. Observe that $s^{k+1}$ is also a point of decay for $\bsign{\Gamma}_{r_k\unit}$. By Lemma \ref{lem_decay_set}, then also the points $\sigma(r)$ with $r \in (r_k,r_{k+1})$ are contained in the decay set of $\bsign{\Gamma}_{r_k\unit}$, and thus in $\Psi(\Gamma)$.%

Property (ii) (coercivity): The piecewise constant functions%
\begin{equation*}
  \tilde{\varphi}_{\min}(r) := \left\{\begin{array}{rl}
	                                       r_k & \mbox{if } r \in [r_k,r_{k+1}) \\
																				   0 & \mbox{if } r = 0
																			\end{array}\right.%
\end{equation*}
and
\begin{equation*}
  \tilde{\varphi}_{\max}(r) := \left\{\begin{array}{rl}
	                                     \varphi(r_{k+1}) & \mbox{if } r \in [r_k,r_{k+1}) \\
																				   0 & \mbox{if } r = 0
																			\end{array}\right.
\end{equation*}
can be used to lower-bound and upper-bound the components of $\sigma$, respectively. By Lemma \ref{lem_mir_A16}, the existence of $\varphi_{\min},\varphi_{\max} \in \KC_{\infty}$ with $\varphi_{\min} \leq \tilde{\varphi}_{\min}$ and $\tilde{\varphi}_{\max} \leq \varphi_{\max}$ follow.%

Property (iii) (monotonicity): This is clear by construction.%

Property (C1): Norm-continuity of $\sigma$ is obvious at any $\rho \notin \{r_k : k \in \Z\}$. If $(\rho_n)_{n\in\N}$ converges to some $r_k$ from below (i.e.~from the left), then $\rho_n \in (r_{k-1},r_k]$ for all sufficiently large $n$. But on $(r_{k-1},r_k]$, $\sigma$ is piecewise linear, hence norm-continuous.%
 
Property (C2): Again, we only have to consider the case $\rho = r_k$ for some $k \in \Z$. If $(\rho_n)_{n\in\N}$ converges to $r_k$ from above (i.e.~from the right), for all sufficiently large $n$ we can find $m_n$ such that $\sigma(r_k) \leq \sigma(\rho_n) \leq \bsign{\Gamma}_{r_k\unit}^{m_n}(s^{k+1})$, and $m_n \rightarrow \infty$ as $n \rightarrow \infty$. Since $\bsign{\Gamma}_{r_k\unit}^{m_n}(s^{k+1}) \stackrel{\star}{\rightharpoonup} \sigma(r_k)$, then also $\sigma(\rho_n) \stackrel{\star}{\rightharpoonup} \sigma(r_k)$.
\end{proof}

If we additionally assume that the fixed point $\sigma^*(r)$ is semi-attractive from above, i.e.~that the convergence $\bsign{\Gamma}_{r\unit}^n(s) \rightarrow \sigma^*(r)$ for all $s \geq \sigma^*(r)$ holds in the norm-topology, then $\sigma$ is obviously continuous.%

We have thus proved the following theorem.%

\begin{theorem}\label{thm_third_path_construction}
Let $\Gamma$ be a generalized gain operator with the following properties:%
\begin{enumerate}
\item[(i)] $\Gamma$ satisfies the $\oplus$-MBI property.%
\item[(ii)] The fixed point $\sigma^*(r)$ of $\bsign{\Gamma}_{r\unit}$ is semi-attractive from above for all $r \geq 0$.%
\end{enumerate}
Then there exists a $C^0$-path of decay for $\Gamma$.
\end{theorem}

\begin{remark}
Observe that assumptions (a) and (b) together imply that $\Sigma(\Gamma)$ is UGAS. Indeed, (a) implies UGS and (b) implies GATT by letting $r = 0$.
\end{remark}

From the construction of $\sigma$, we obtain a new characterization of the $\oplus$-MBI property.%

\begin{theorem}\label{thm_maxmbi_path_char}
For a generalized gain operator $\Gamma$, the following statements are equivalent:%
\begin{enumerate}
\item[(a)] There exists $\rho \in \KC_{\infty}$ such that $\Gamma_{\rho}$ satisfies the $\oplus$-MBI property.%
\item[(b)] There exist $\rho \in \KC_{\infty}$ and a path $\sigma:\R_+ \rightarrow \ell^{\infty}_+(\IC)$ satisfying the properties listed in Proposition \ref{prop_maxmbi_path} for the operator $\Gamma_{\rho}$ (in place of $\Gamma$).
\end{enumerate}
\end{theorem}

\begin{proof}
We already know that (a) implies (b). To prove the converse, we can essentially use the same arguments as in the proof of Theorem \ref{thm_necessary_conditions}. Let $s > 0$ (without loss of generality) and $r \geq 0$ be chosen such that $s \leq r\unit \oplus \Gamma_{\rho}(s)$ (for a $\rho$ smaller than that in (b)). We define%
\begin{equation*}
  r' := \min\{ \tilde{r} \geq 0 : s \leq \sigma(\tilde{r}) \}.%
\end{equation*}
To show that $r'$ is well-defined, we prove that the set $A := \{\tilde{r} \geq 0 : s \leq \sigma(\tilde{r})\}$ is closed in $\R_+$. Hence, let $(\tilde{r}_n)_{n\in\N}$ be a sequence in $A$, converging to some $\tilde{r} \in \R_+$. We have $s \leq \sigma(\tilde{r}_n)$ for all $n$. From the continuity assumptions on $\sigma$, it follows that $\sigma$ is sequentially continuous in the weak$^*$-topology. Hence, $\sigma(\tilde{r}_n) \stackrel{\star}{\rightharpoonup} \sigma(\tilde{r})$, which clearly implies $s \leq \sigma(\tilde{r})$. It follows that $\tilde{r} \in A$, showing that $A$ is closed. Moreover, we must have $r' > 0$, since otherwise $s = 0$. Let us assume to the contrary that there exists an $\ep > 0$ with $s \leq \sigma(r') - \ep\unit$. In this case, let $(r_n)_{n\in\N}$ be a sequence in $[0,r')$, converging to $r'$. Since $\sigma$ is norm-continuous from below, it follows that $\|\sigma(r_n) - \sigma(r')\| \leq \ep$ for $n$ sufficiently large. This implies $s \leq \sigma(r') - \ep\unit \leq \sigma(r_n)$, which contradicts the minimality of $r'$. Hence, for every $\ep > 0$ there exists an index $i \in \IC$ with $\sigma_i(r') - s_i < \ep$. The rest of the proof is identical to that of statement (b) in Theorem \ref{thm_necessary_conditions}.
\end{proof}

Let us now assume that $\IC$ is infinite and $\JC \subset \IC$ is a nonempty finite subset. This subset corresponds to a finite subnetwork with associated gain operator%
\begin{equation*}
  \Gamma^{\langle\JC\rangle}:\ell^{\infty}_+(\JC) \rightarrow \ell^{\infty}_+(\JC),\quad \Gamma^{\langle \JC \rangle}_i(s) := \Gamma_i(s_{|\IC_i \cap \JC}).%
\end{equation*}
The associated graph $\GC^{\JC}$ has vertex set $V(\GC^{\JC}) = \JC$ and edge set $E(\GC^{\JC}) = \{ ji : i \in \JC,\ j \in \IC_i \cap \JC \}$.%

Now, assume that $\Gamma$ admits a path $\sigma:\R_+ \rightarrow \ell^{\infty}_+(\IC)$ satisfying all properties listed in Proposition \ref{prop_maxmbi_path}. Then, we can consider the restricted path%
\begin{equation*}
  \sigma^{\JC}(r) := (\sigma_i(r))_{i \in \JC},\quad \sigma^{\JC}:\R_+ \rightarrow \ell^{\infty}_+(\JC).%
\end{equation*}
Since $\JC$ is finite, this path is continuous, and moreover, it is increasing and inherits the coercivity estimates of $\sigma$. Finally, for every $i \in \JC$ we have%
\begin{equation*}
  \Gamma^{\langle\JC\rangle}_i(\sigma^{\JC}(r)) = \Gamma_i( \sigma(r)_{|\IC_i \cap \JC} ) \leq \Gamma_i(\sigma(r)) \leq (\id + \rho)^{-1}(\sigma_i(r)).%
\end{equation*}
This implies $\Gamma^{\langle\JC\rangle}_{\rho} \circ \sigma^{\JC} \leq \sigma^{\JC}$, hence $\sigma^{\JC}$ is a path of strict decay for $\Gamma^{\langle \JC \rangle}$. We thus have the following result.%

\begin{corollary}\label{cor_finite_subnetworks}
Assume that $\Gamma_{\rho}$ satisfies the $\oplus$-MBI property for some $\rho \in \KC_{\infty}$. Then the gain operator of every finite subnetwork admits a path of strict decay and the coercivity bounds are the same for all subnetworks.
\end{corollary}

\section{Main results}\label{sec_final_theorem}

Putting our partial results together allows us to formulate the following theorem.%

\begin{theorem}\label{thm_final}
For a generalized gain operator $\Gamma$, consider the following statements:%
\begin{enumerate}
\item[(a)] There exists a path of strict decay for $\Gamma$.%
\item[(b)] There exists a $C^0$-path of strict decay for $\Gamma$.%
\item[(c)] There exists $\rho \in \KC_{\infty}$ such that $\Sigma(\Gamma_{\rho})$ is UGAS and $\Gamma_{\rho}$ satisfies the $\oplus$-MBI property.%
\item[(d)] There exists $\rho \in \KC_{\infty}$ such that $\Gamma_{\rho}$ satisfies both uniform NJI conditions.%
\item[(e)] There exists a path $\sigma:\R_+ \rightarrow \ell^{\infty}_+(\IC)$, satisfying all properties of a $C^0$-path of strict decay except for norm-continuity from above, which is replaced by sequential weak$^*$ continuity (from above).%
\item[(f)] For every finite subnetwork, there exists a path of strict decay and the coercivity bounds are the same for all subnetworks.%
\end{enumerate}
Then, (a) $\Leftrightarrow$ (b) $\Rightarrow$ (c) $\Leftrightarrow$ (d) $\Rightarrow$ (e) $\Rightarrow$ (f). For the implication (d) $\Rightarrow$ (c), the additional assumption that the operator is locally finite and $|\IC_i|$ is uniformly bounded is needed.
\end{theorem}

\begin{proof}
(a) $\Rightarrow$ (b): Obviously, every path of strict decay is a $C^0$-path of strict decay.%

(b) $\Rightarrow$ (a): This is Proposition \ref{prop_path_of_decay_2}.%

(b) $\Rightarrow$ (c): This is Theorem \ref{thm_necessary_conditions}.%

(c) $\Rightarrow$ (d): This follows from Proposition \ref{prop_GATT_implies_uniform_NJI} and Proposition \ref{prop_maxmbi_nji_char}.%

(d) $\Rightarrow$ (c): By Proposition \ref{prop_maxmbi_nji_char}, the SU-NJI condition is equivalent to the $\oplus$-MBI property. The $\oplus$-MBI property implies that the corresponding decay set is cofinal by Proposition \ref{prop_maxmbi_char}. This together with the TU-NJI condition implies that $\Sigma(\Gamma_{\rho})$ is GATT for some $\rho$ by Proposition \ref{prop_uniform_NJI_implies_GATT}. The $\oplus$-MBI property also implies that $\Sigma(\hat{\Gamma}_{\rho})$, and thus also $\Sigma(\Gamma_{\rho})$, is UGS, again by Proposition \ref{prop_maxmbi_char}. UGS and GATT together are equivalent to UGAS, which completes the proof.%

(d) $\Rightarrow$ (e): This follows from Theorem \ref{thm_maxmbi_path_char}. Note that here only the SU-NJI condition is needed as an assumption.%

(e) $\Rightarrow$ (f): This is Theorem \ref{thm_maxmbi_path_char} together with Corollary \ref{cor_finite_subnetworks}.
\end{proof}

The remaining question is whether we can get back from (c) to (b), or which additional assumptions are needed for this implication. The hope would be that statements (a)--(d) are equivalent.%

Some of the results that will be shown in the next section show that the implication ``UGAS $\Rightarrow$ $\oplus$-MBI property'' sometimes holds under additional assumptions on the gain operator. For instance, this is true if $\Gamma$ is uniformly continuous and $|\IC| < \infty$. In general, however, we expect that the implication does not hold.

Using Proposition \ref{prop_maxmbi_char}, Theorem \ref{thm_maxmbi_path_char} and Proposition \ref{prop_maxmbi_nji_char}, we obtain the following weaker version of Theorem \ref{thm_final}.%

\begin{theorem}\label{thm_final_weak}
For a generalized gain operator $\Gamma$, the following statements are equivalent:%
\begin{enumerate}
\item[(a)] There exists a path $\sigma:\R_+ \rightarrow \ell^{\infty}_+(\IC)$, satisfying all properties of a $C^0$-path of strict decay except for norm-continuity from above, which is replaced by sequential weak$^*$ continuity (from above).%
\item[(b)] There exists $\rho \in \KC_{\infty}$ such that $\Gamma_{\rho}$ satisfies the $\oplus$-MBI property.%
\item[(c)] There exists $\rho \in \KC_{\infty}$ such that $\Gamma_{\rho}$ satisfies the SU-NJI condition.%
\item[(d)] There exists $\rho \in \KC_{\infty}$ such that $\Sigma(\hat{\Gamma}_{\rho})$ is UGS, $\Sigma(\Gamma_{\rho})$ is GATT$^{\,*}$ and the set $\bigcup_{r \geq 0}\Fix(\bsign{\Gamma}_{\rho,r\unit})$ is coercive.%
\end{enumerate}
\end{theorem}

This theorem becomes particularly interesting in the finite case, where it can be formulated as follows. Recall that here $\Gamma$ is simply a continuous monotone operator with $\Gamma(0) = 0$.%

\begin{theorem}\label{thm_final_weak_finite_case}
If $|\IC| < \infty$, then the following statements are equivalent:%
\begin{enumerate}
\item[(a)] There exists a path of strict decay for $\Gamma$.%
\item[(b)] There exists $\rho \in \KC_{\infty}$ such that $\Gamma_{\rho}$ satisfies the $\oplus$-MBI property.%
\item[(c)] There exists $\rho \in \KC_{\infty}$ such that $\Gamma_{\rho}$ satisfies the SU-NJI condition.%
\item[(d)] There exists $\rho \in \KC_{\infty}$ such that $\Sigma(\Gamma_{\rho})$ is UGAS, $\Sigma(\hat{\Gamma}_{\rho})$ is UGS and the set $\bigcup_{r \geq 0}\Fix(\bsign{\Gamma}_{\rho,r\unit})$ is coercive.%
\end{enumerate}
\end{theorem}

In Subsection \ref{subsec_finite_case}, we will show that for (non-generalized) gain operators the SU-NJI condition can be reduced to the usual NJI condition under an additional assumption.

\section{Particular classes of gain operators}\label{sec_classes}

In this section, we take a closer look at four particular classes of gain operators.%

\subsection{Uniformly continuous operators}\label{subsec_uniformly_continuous}

By definition, any gain operator is uniformly continuous on bounded sets. In this subsection, we will show that the assumption of uniform continuity on the whole domain of definition allows for stronger conclusions. In particular, under this assumption we will obtain more clarity about the relation between the two necessary conditions for the existence of a path of strict decay, namely UGAS of $\Sigma(\Gamma_{\rho})$ and the $\oplus$-MBI property of $\Gamma_{\rho}$, or the corresponding NJI conditions.%

\begin{lemma}\label{lem_unifcont_op}
Given a uniformly continuous gain operator $\Gamma$, for each $n \in \N$ there exists a $\varphi \in \KC_{\infty}$ such that for all $\ep > 0$ and $s \in \ell^{\infty}_+(\IC)$, $i \in \IC$ the following holds: $\Gamma_j(s) \geq s_j - \varphi(\ep)$ for all $j \in \NC^-_i(n-1)$ implies $\Gamma_i^n(s) \geq s_i - \ep$.
\end{lemma}

\begin{proof}
We use the characterization of uniform continuity in terms of comparison functions. Let $\gamma \in \KC_{\infty}$ be chosen such that%
\begin{equation*}
  \|\Gamma(s^1) - \Gamma(s^2)\| \leq \gamma(\|s^1 - s^2\|)%
\end{equation*}
holds for all $s^1,s^2 \in \ell^{\infty}_+(\IC)$. Given $n \in \N$ and $\ep > 0$, we then choose $\ep_i,\delta_i > 0$, $1 \leq i \leq n - 1$ as in the proof of Lemma \ref{lem_aux_stability}, now using the function $\gamma$. That is, we put $\delta_n := \ep$, $\ep_{n-1} := \delta_n/2$ and continue recursively by%
\begin{equation*}
  \delta_l := \gamma^{-1}(\ep_l), \quad \ep_{l-1} := \frac{\delta_l}{2}%
\end{equation*}
We do this for $l = n - 1$ down to $l = 2$. Eventually, we choose $\delta_1 =: \delta$ as%
\begin{equation*}
  \delta_1 := \min\left\{ \gamma^{-1}(\ep_1), \ep_1,\ep_2,\ldots,\ep_{n-1} \right\}.%
\end{equation*}
Without loss of generality, we can actually assume that $\gamma^{-1} < \id$, which implies $\delta_1 = \gamma^{-1}(\ep_1)$. This obviously defines a $\KC_{\infty}$-function $\varphi$, only depending on $n$, such that $\delta = \varphi(\ep)$. The rest of the proof is completely analogous to that of Lemma \ref{lem_aux_stability}.
\end{proof}

What makes Lemma \ref{lem_unifcont_op} superior to Lemma \ref{lem_aux_stability} is that its statement is not only meaningful for $\ep$ close to zero. The following proposition shows how this can be used.%

\begin{proposition}\label{prop_unifcont_gainop}
Let $\Gamma$ be a uniformly continuous gain operator. Then, under each of the following assumptions, $\Gamma$ satisfies the $\oplus$-MBI property:%
\begin{enumerate}
\item[(i)] $|\IC| < \infty$ and $\Sigma(\Gamma)$ is GATT.%
\item[(ii)] $\Sigma(\Gamma)$ is UGAS with a $\KC\LC$-function of the form $\beta(r,n) = r\gamma(n)$, $\gamma \in \LC$.%
\end{enumerate}
\end{proposition}

\begin{proof}
(i) We prove that $\Gamma$ satisfies the SU-NJI condition with a $\KC_{\infty}$-function $\varphi$. To this end, we pick $m \in \N$ such that $\NC^-_i = \NC^-_i(n)$ for all $i \in \IC$ and $n \geq m$, which is possible due to the finiteness of $\IC$. Then, using Lemma \ref{lem_unifcont_op}, we choose $\tilde{\varphi} = \tilde{\varphi}(m) \in \KC_{\infty}$ accordingly. We put $\varphi(r) :\equiv \tilde{\varphi}(r/2)$ and assume towards a contradiction that there exist $s \geq 0$ and $i \in \IC$ with $s_i > 0$, but $s_j < \varphi(s_i)$ or $\Gamma_j(s) \geq s_j$ for all $j \in \NC^-_i$. Observe that in both cases $\Gamma_j(s) \geq s_j - \tilde{\varphi}(s_i/2)$. Hence, Lemma \ref{lem_unifcont_op} implies $\Gamma_i^n(s) \geq s_i/2 > 0$ for all $n \geq m$. This obviously contradicts GATT.%

(ii) We fix an $m \in \N$ with $\gamma(m) < 1/4$, choose $\tilde{\varphi} = \tilde{\varphi}(m) \in \KC_{\infty}$ according to Lemma \ref{lem_unifcont_op}, and put $\varphi(r) :\equiv \tilde{\varphi}(r/2)$. We show the following strong version of the TU-NJI condition: for all $s \geq 0$ and $i \in \IC$:%
\begin{equation*}
  s_i \geq \frac{\|s\|}{2} \quad \Rightarrow \quad \exists j \in \NC^-_i(m):\ s_j \geq \varphi\left(\frac{\|s\|}{2}\right) \wedge \Gamma_j(s) < s_j.%
\end{equation*}
Assume towards a contradiction that there are $s \geq 0$ and $i \in \IC$ with $s_i \geq \|s\|/2$, but $s_j < \varphi(\|s\|/2)$ or $\Gamma_j(s) \geq s_j$ for all $j \in \NC^-_i(m)$; in both cases, $\Gamma_j(s) \geq s_j - \varphi(\|s\|/2)$ holds. Hence, Lemma \ref{lem_unifcont_op} implies $\Gamma_i^m(s) \geq s_i - \|s\|/4$. Then%
\begin{equation*}
  \frac{\|s\|}{2} \leq s_i \leq \Gamma_i^m(s) + \frac{\|s\|}{4} \leq \|\Gamma^m(s)\| + \frac{\|s\|}{4} \leq \|s\|\gamma(m) + \frac{\|s\|}{4} < \frac{\|s\|}{2},%
\end{equation*}
which is a contradiction. Now, towards another contradiction, assume that $s \leq r\unit \oplus \Gamma(s)$ and $\|s\| > 2\varphi^{-1}(r)$. Pick any $i \in \IC$ with $s_i \geq \|s\|/2$. Then, by the above, there exists $j \in \NC^-_i(m)$ with $s_j \geq \varphi(\|s\|/2)$ and $\Gamma_j(s) < s_j$. The latter implies $s_j \leq r$. Altogether, $r < \varphi(\|s\|/2) \leq s_j \leq r$, which is a contradiction. Hence, we have shown that $s \leq r\unit \oplus \Gamma(s)$ implies $\|s\| \leq 2\varphi^{-1}(r)$.
\end{proof}

\subsection{Max-type operators}\label{subsec_maxtype}

Recall that a gain operator $\Gamma$ is called a \emph{max-type operator} if all MAFs $\mu_i$ are of the form $\mu_i(s) = \sup_{j\in\IC}s_j = \|s\|$. In this case, $\Gamma$ satisfies%
\begin{equation*}
  \Gamma(s^1 \oplus s^2) = \Gamma(s^1) \oplus \Gamma(s^2) \mbox{\quad for all\ } s^1,s^2 \in \ell^{\infty}_+(\IC).%
\end{equation*}
We also call a generalized gain operator a max-type gain operator if it has this property. Gain operators of this type are particularly well-behaved. Using our general results, it is easy to provide a full characterization of the existence of a path of strict decay.%

\begin{theorem}\label{thm_maxtype_main}
For a max-type gain operator $\Gamma$, the following statements are equivalent:%
\begin{enumerate}
\item[(a)] There exists a path of strict decay for $\Gamma$.%
\item[(b)] There exists $\rho \in \KC_{\infty}$ such that $\Sigma(\Gamma_{\rho})$ is UGAS.%
\end{enumerate}
If there exists $M > 0$ with $|\IC_i| \leq M$ for all $i \in \IC$, then another equivalent statement is the following:%
\begin{enumerate}
\item[(c)] There exists $\rho \in \KC_{\infty}$ such that $\Gamma_{\rho}$ satisfies both uniform NJI conditions.
\end{enumerate}
\end{theorem}

\begin{proof}
Recall that the implications ``(a) $\Rightarrow$ (c)'', ``(a) $\Rightarrow$ (b)'' and ``(c) $\Rightarrow$ (b)'' also hold in the general case. It remains to prove that (b) implies (a). We first prove that UGAS of $\Sigma(\Gamma_{\rho})$ implies that $\Gamma_{\rho}$ satisfies the $\oplus$-MBI property. Hence, let $s \leq b \oplus \Gamma_{\rho}(s) = \bsign{\Gamma}_{\rho,b}(s)$ for some $b,s \in \ell^{\infty}_+(\IC)$. This implies $s \leq \bsign{\Gamma}_{\rho,b}^n(s)$ for all $n \in \N$. Observe that $\Gamma_{\rho}$ is also a max-type operator, which easily implies%
\begin{equation}\label{eq_maxtype_iterates}
  \bsign{\Gamma}_{\rho,b}^n(s) \equiv \Gamma_{\rho}^n(s) \oplus \bigoplus_{k=0}^{n-1} \Gamma^k_{\rho}(b) \mbox{\quad for all\ } n \geq 0.%
\end{equation}
By assumption, $\Gamma_{\rho}^n(s) \rightarrow 0$. Hence, $s \leq \bigoplus_{k=0}^{\infty} \Gamma^k_{\rho}(b)$, leading to%
\begin{equation*}
  \|s\| \leq \sup_{k \in \N_0} \|\Gamma^k_{\rho}(b)\| \leq \sup_{k\in\N_0} \beta(\|b\|,k) = \beta(\|b\|,0) =: \varphi(\|b\|).%
\end{equation*}
Thus, $\Gamma_{\rho}$ satisfies the $\oplus$-MBI property. To obtain the desired path, we can now try to apply Theorem \ref{thm_third_path_construction}, showing that the fixed point $\sigma^*(r)$ of $\bsign{\Gamma}_{\rho,r\unit}$ is semi-attractive from above. We obtain even more. Putting $b = r\unit$ in \eqref{eq_maxtype_iterates} implies that $\bsign{\Gamma}_{\rho,r\unit}$ only has one fixed point, namely $\bigoplus_{k=0}^{\infty} \Gamma^k_{\rho}(b)$, which is globally attractive, since $\Gamma_{\rho}^n(s)$ converges to zero. 
\end{proof}

\begin{remark}
The theorem improves \cite[Thm.~VI.1]{KMZ}, where uniform Lipschitz conditions on $\gamma_{ij}$ and $\gamma_{ij}^{-1}$ were required to obtain a path of strict decay. Moreover, \cite[Thm.~VI.1]{KMZ} only provided a sufficient condition and no relation to an NJI condition.
\end{remark}

\subsection{Homogeneous, subadditive operators}\label{subsec_homogeneous}

We call a generalized gain operator $\Gamma$ \emph{homogeneous (of degree one)} if $\Gamma(\alpha s) = \alpha \Gamma(s)$ for all $s \in \ell^{\infty}_+(\IC)$ and $\alpha \geq 0$. We call $\Gamma$ \emph{subadditive} if $\Gamma(s^1 + s^2) \leq \Gamma(s^1) + \Gamma(s^2)$ for all $s^1,s^2 \in \ell^{\infty}_+(\IC)$. In this subsection, we consider operators which are both homogeneous and subadditive. Examples are sum-type and max-type operators with linear interconnection gains. In the first case, $\Gamma$ is the restriction of a linear operator. In the second case, we speak of a max-linear operator. Except for the equivalence to statement (d), the next theorem is known, see \cite[Prop.~8]{Mea}, but we give an alternative proof.%

\begin{theorem}\label{thm_homogeneous_subadditive}
For a homogeneous and subadditive gain operator $\Gamma$, the following statements are equivalent:%
\begin{enumerate}
\item[(a)] There exists a path of strict decay for $\Gamma$.%
\item[(b)] The system $\Sigma(\Gamma)$ is UGAS.%
\item[(c)] The following ``spectral radius'' condition holds:%
\begin{equation}\label{eq_spectralradius_condition}
  \inf_{n \in \N} \|\Gamma^n(\unit)\| < 1.%
\end{equation}
\end{enumerate}
If there exists $M > 0$ with $|\IC_i| \leq M$ for all $i \in \IC$, then another equivalent statement is the following:%
\begin{enumerate}
\item[(d)] There exists $\rho \in \KC_{\infty}$ such that $\Gamma_{\rho}$ satisfies both uniform NJI conditions.
\end{enumerate}
\end{theorem}

\begin{proof}
The implications ``(a) $\Rightarrow$ (b)'', ``(d) $\Rightarrow$ (b)'' and ``(a) $\Rightarrow$ (d)'' also hold in the general case. The equivalence ``(b) $\Leftrightarrow$ (c)'' was shown in \cite[Prop.~8]{Mea}, and we do not repeat the proof here. It remains to (b) implies (a).%

First, we refer to \cite[Prop.~9]{Mea}, which shows that (b) is equivalent to uniform global exponential stability (UGES), which is UGAS with a $\KC\LC$-function of the form $\beta(r,n) = Ma^n r$, where $M > 0$ and $a \in (0,1)$. From the homogeneity of $\Gamma$, it follows that we can choose a linear $\KC_{\infty}$-function $\rho$ (now regarded as a positive real number) such that $\Gamma_{\rho} = (1 + \rho) \cdot \Gamma$ is still UGES. In particular, this implies that (b) holds with $\Gamma_{\rho}$ in place of $\Gamma$. We can apply Proposition \ref{prop_unifcont_gainop} to show that $\Gamma_{\rho}$ satisfies the $\oplus$-MBI property. To this end, we have to show that $\Gamma$ is uniformly continuous. First note that $0 \leq s^1 \leq s^2$ implies%
\begin{equation*}
  \Gamma(s^2) - \Gamma(s^1) = \Gamma(s^1 + (s^2 - s^1)) - \Gamma(s^1) \leq \Gamma(s^2 - s^1)%
\end{equation*}
and thus, by continuity at $s = 0$, for every $\ep > 0$ we find $\delta > 0$ such that $\|s^2 - s^1\| \leq \delta$ implies%
\begin{equation*}
  \|\Gamma(s^2) - \Gamma(s^1)\| \leq \|\Gamma(s^2 - s^1)\| = \|\Gamma(s^2 - s^1) - \Gamma(0)\| \leq \ep.%
\end{equation*}
For general $s^1,s^2 \in \ell^{\infty}_+(\IC)$ (not necessarily ordered), we have%
\begin{equation*}
  \Gamma(s^2) - \Gamma(s^1) \leq \Gamma(s^1 \oplus s^2) - \Gamma(s^1) \leq \Gamma(s^1 \oplus s^2 - s^1)%
\end{equation*}
which implies%
\begin{equation*}
  \|\Gamma(s^2) - \Gamma(s^1)\| \leq \|\Gamma(s^1 \oplus s^2 - s^1) - \Gamma(0)\|.%
\end{equation*}
Observing that $\|s^1 \oplus s^2 - s^1\| \leq \|s^2 - s^1\|$, this concludes the proof of uniform continuity. Hence, $\Gamma_{\rho}$ satisfies the $\oplus$-MBI property.%

We claim that for any $r \geq 0$, points $s^1 \leq s^2$ and $n \in \N_0$, we have%
\begin{equation}\label{eq_subadd_contraction}
  \bsign{\Gamma}^n_{r\unit}(s^2) - \bsign{\Gamma}^n_{r\unit}(s^1) \leq \Gamma^n(s^2 - s^1).%
\end{equation}
For $n = 0$, this trivially holds. Assuming that it holds for a fixed $n$, with Lemma \ref{lem_maxdiff} we obtain%
\begin{align*}
  \bsign{\Gamma}^{n+1}_{r\unit}(s^2) - \bsign{\Gamma}^{n+1}_{r\unit}(s^1) &= r\unit \oplus \Gamma(\bsign{\Gamma}_{r\unit}^n(s^2)) - r\unit \oplus \Gamma(\bsign{\Gamma}_{r\unit}^n(s^1)) \\
	&\leq \Gamma(\bsign{\Gamma}_{r\unit}^n(s^2)) - \Gamma(\bsign{\Gamma}_{r\unit}^n(s^1)) \leq \Gamma(\bsign{\Gamma}_{r\unit}^n(s^2) - \bsign{\Gamma}_{r\unit}^n(s^1)) \\
	&\leq \Gamma(\Gamma^n(s^2 - s^1)) = \Gamma^{n+1}(s^2 - s^1).%
\end{align*}
Here, the second inequality follows from subadditivity of $\Gamma$ and the last one from the induction hypothesis. The same estimates hold with $\Gamma_{\rho}$ in place of $\Gamma$, in case that $\rho$ is linear. Hence,%
\begin{equation*}
  \|\bsign{\Gamma}^n_{\rho,r\unit}(s^2) - \bsign{\Gamma}^n_{\rho,r\unit}(s^1)\| \leq Ma^n \|s^2 - s^1\|%
\end{equation*}
for appropriate constants $M > 0$ and $a \in (0,1)$. This implies that $\bsign{\Gamma}_{\rho,r\unit}$ can only have one fixed point and every trajectory converges to this fixed point at exponential speed. Hence, Theorem \ref{thm_third_path_construction} yields a path of strict decay for $\Gamma$.
\end{proof}

We finally prove the following theorem, which is somewhat surprising, since it shows that the $\oplus$-MBI property is sufficient to obtain a path of strict decay for homogeneous subadditive operators.%

\begin{theorem}\label{thm_homogenenous_subadditive_2}
Let $\Gamma$ be a homogeneous and subadditive gain operator and assume that $\Gamma_{\rho}$ satisfies the SU-NJI condition with a linear $\rho \in \KC_{\infty}$. Then there exists a linear path of strict decay for $\Gamma$.
\end{theorem}

\begin{proof}
First, observe that $\Gamma_{\rho}$ is a subadditive and homogeneous gain operator if $\rho$ is linear. Consider the inequality \eqref{eq_subadd_contraction}. Our assumption (which is equivalent to the $\oplus$-MBI property), implies that $\Gamma_{\rho}^n(s^2 - s^1)$ converges to the origin componentwise, enforcing the uniqueness of the fixed point of $\bsign{\Gamma}_{\rho,r\unit}$. We consider the path $\sigma_*$, which maps each $r$ to the unique fixed point of $\bsign{\Gamma}_{\rho,r\unit}$. Then, we have the identity%
\begin{equation*}
  \sigma_*(r) \equiv r\unit \oplus \Gamma_{\rho}(\sigma_*(r)).%
\end{equation*}
Multiplying both sides by a nonnegative scalar $\alpha$ yields%
\begin{equation*}
  \alpha\sigma_*(r) \equiv (\alpha r)\unit \oplus \Gamma_{\rho}(\alpha\sigma_*(r))%
\end{equation*}
due to homogeneity. By uniqueness of the fixed point, this implies $\sigma_*(\alpha r) = \alpha \sigma_*(r)$. We thus have $\sigma_*(r) \equiv r \sigma_*(1)$, showing that $\sigma_*$ is linear, and hence continuous. Obviously, $\sigma_*$ also satisfies the other three properties of a path of strict decay (see the proof of Proposition \ref{prop_hatsigma_ugs}).
\end{proof}

\begin{remark}
We note that several further characterizations of the equivalent conditions in Theorem \ref{thm_homogeneous_subadditive} are available, especially for linear operators, see \cite{GMi}. This includes, in particular, the spectral radius condition $r(\Gamma) < 1$, which is an immediate consequence of \eqref{eq_spectralradius_condition} and Gelfand's formula.
\end{remark}

\subsection{Finite-dimensional operators}\label{subsec_finite_case}

In this subsection, we study the finite case, where $\ell^{\infty}(\IC)$ can be identified with the Euclidean space $\R^N$, $N := |\IC|$, and $\ell^{\infty}_+(\IC)$ with $\R^N_+ := \{ s \in \R^N : s_i \geq 0,\ 1 \leq i \leq N\}$. In particular, we clarify how our results are related to previous results from the literature.%

First, let us clarify which assumptions we have to impose on the functions $\gamma_{ij}$ and $\mu_i$ to obtain a gain operator in the sense of Definition \ref{def_normal_gain_op}.%

\begin{proposition}\label{prop_finitedim_gainop}
Let $\GC = (V(\GC),E(\GC))$ be a finite directed graph without self-loops (i.e.~edges from a vertex to itself) and put $\IC := V(\GC)$. Consider families $\{\gamma_{ij} : ji \in E(\GC)\}$ and $\{\mu_i : 1 \leq i \leq |\IC|\}$ with $\gamma_{ij} \in \KC_{\infty}$ and $\mu_i:\ell^{\infty}_+(\IC) \rightarrow \R_+$ satisfying the following properties for all $i \in \IC$:%
\begin{enumerate}
\item[(i)] $\mu_i(0) = 0$ and $\mu_i(s) > 0$ if $s > 0$.%
\item[(ii)] $0 \leq s^1 \leq s^2$ implies $\mu_i(s^1) \leq \mu_i(s^2)$.%
\item[(iii)] $\|s\| \rightarrow \infty$ implies $\mu_i(s) \rightarrow \infty$.%
\item[(iv)] $\mu_i$ is continuous.%
\end{enumerate}
Then the operator $\Gamma:\ell^{\infty}_+(\IC) \rightarrow \ell^{\infty}_+(\IC)$, defined by $\Gamma_i(s) := \mu_i([\gamma_{ij}(s_j)]_{j\in\IC_i})$ for all $i \in \IC$, is a gain operator.
\end{proposition}

\begin{proof}
Since $\{\gamma_{ij}\}$ is a finite family of functions, it is pointwise equicontinuous. It remains to verify (M1)--(M4) for $\{\mu_i\}$. This, in turn, reduces to showing that $\mu_i(s) \geq \xi(\|s\|)$ for some $\xi \in \KC_{\infty}$, since the rest trivially follows from our assumptions. To this end, we define%
\begin{equation*}
  \tilde{\xi}(r) := \min_{(i,j) \in \IC \tm \IC}\mu_i(re^j).%
\end{equation*}
From assumption (i), we obtain $\tilde{\xi}(0) = 0$. Clearly, $\tilde{\xi}(r) \geq 0$ for all $r \in \R_+$. From assumption (iv), it follows, moreover, that $\tilde{\xi}$ is continuous. If $r_1 \leq r_2$, then $r_1 e^j \leq r_2 e^j$ for each $j$, and hence assumption (ii) implies $\tilde{\xi}(r_1) \leq \tilde{\xi}(r_2)$, showing that $\tilde{\xi}$ is non-decreasing. If $r \rightarrow \infty$, then $\|re^j\| \rightarrow \infty$, and hence, following from assumption (iii), $\tilde{\xi}(r) \rightarrow \infty$. Finally, assumption (i) implies $\tilde{\xi}(r) > 0$ if $r > 0$. By Lemma \ref{lem_mir_A16}, we can then bound $\tilde{\xi}$ from below by a $\KC_{\infty}$-function $\xi$. To show that $\mu_i(s) \geq \xi(\|s\|)$, pick $s \in \ell^{\infty}_+(\IC)$ and choose an index $j$ such that $s_j = \|s\|$. Then $s \geq s_j e^j = \|s\|e^j$, implying $ \mu_i(s) \geq \mu_i(\|s\|e^j) \geq \tilde{\xi}(\|s\|) \geq \xi(\|s\|)$. This concludes the proof.
\end{proof}

\begin{remark}
The assumptions of the proposition can also be found in the literature on finite-dimensional gain operators.
\end{remark}

For our general result, we need an additional \emph{uniform continuity} assumption for the MAFs that we have not used before. We can formulate this assumption in terms of a comparison function as follows: There exists a $\KC_{\infty}$-function $\gamma$ such that $0 \leq s^1 \leq s^2$ implies%
\begin{equation}\label{eq_mafs_uniform_cont}
  \mu_i(s^2) - \mu_i(s^1) \leq \gamma(\|s^2 - s^1\|) \mbox{\quad for all\ } i \in \IC.
\end{equation}
Observe that this is indeed equivalent to the uniform continuity of each $\mu_i$. First, since there are only finitely many MAFs, the same $\gamma$ can be used for all $i \in \IC$. Second, if $s^1$ and $s^2$ are arbitrary elements of $\ell^{\infty}_+(\IC)$ (not necessarily ordered), it follows from the assumption that%
\begin{equation*}
  \mu_i(s^2) - \mu_i(s^1) \leq \mu_i(s^1 \oplus s^2) - \mu_i(s^1) \leq \gamma(\|s^1 \oplus s^2 - s^1\|) \leq \gamma(\|s^2 - s^1\|)%
\end{equation*}
and analogously for $\mu_i(s^1) - \mu_i(s^2)$. Hence, $|\mu_i(s^2) - \mu_i(s^1)| \leq \gamma(\|s^2 - s^1\|)$, which is the characterization of uniform continuity in terms of comparison functions.%

In the literature on finite-dimensional small-gain theory (see \cite{DR1,DR2,Ruf}), instead the assumption of \emph{subadditivity} is used, which is a special case of the above. To see this, note that a MAF $\mu$ is subadditive if $\mu(s^1 + s^2) \leq \mu(s^1) + \mu(s^2)$ for all $s^1,s^2 \in \ell^{\infty}_+(\IC)$. If $0 \leq s^1 \leq s^2$, this implies%
\begin{equation*}
  \mu(s^2) - \mu(s^1) = \mu(s^1 + (s^2 - s^1)) - \mu(s^1) \leq \mu(s^2 - s^1) \leq \mu( \|s^2 - s^1\|\unit )%
\end{equation*}
and we can define $\gamma(r) := \mu(r\unit)$, which can be bounded from above by a $\KC_{\infty}$-function.%

\begin{theorem}\label{thm_finite_general_case}
Assume that $|\IC| < \infty$ and consider a gain operator with uniformly continuous MAFs. Then the following statements are equivalent:%
\begin{enumerate}
\item[(a)] There exists a path of strict decay for $\Gamma$.%
\item[(b)] There exists $\rho \in \KC_{\infty}$ such that $\Sigma(\Gamma_{\rho})$ is UGAS.
\item[(c)] There exists $\rho \in \KC_{\infty}$ such that $\Gamma_{\rho}$ satisfies the $\oplus$-MBI property.%
\item[(d)] There exists $\rho \in \KC_{\infty}$ such that $\Gamma_{\rho}$ satisfies the NJI condition.%
\end{enumerate}
\end{theorem}

\begin{proof}
Looking at our general results in Section \ref{sec_final_theorem} and observing that (b) trivially implies (d), we can easily see that the proof is completed by showing that (d) implies (c). To this end, consider $s \in \ell^{\infty}_+(\IC)$ and $r \geq 0$ with $s \leq r\unit \oplus \Gamma(s)$. Without loss of generality, we can assume that $s > 0$. We split the index set $\IC$ into the two subsets%
\begin{equation*}
  \IC^1 := \left\{ i \in \IC : s_i \leq r \right\},\quad \IC^2 := \IC \backslash \IC^1.%
\end{equation*}
By assumption, there exists at least one index $i$ with $\Gamma_i(s) < s_i$, which implies $s_i \leq r$. Hence, $\IC^1$ is nonempty. We define the $\KC_{\infty}$-function%
\begin{equation*}
  \varphi_1 := (\rho \circ \xi \circ \eta)^{-1} \circ \gamma \circ \left(\max_{ji \in E(\GC)}\gamma_{ij}\right),%
\end{equation*}
where $\rho$ comes from (a), $\xi$ comes from (M1), $\eta \leq \gamma_{ij}$ for all $ji \in E(\GC)$ with $\eta \in \KC_{\infty}$ (such $\eta$ trivially exists), and $\gamma$ comes from \eqref{eq_mafs_uniform_cont}. Without loss of generality, we can assume that%
\begin{equation}\label{eq_finite_general_wlogs}
  \varphi_1(r) > r \mbox{\quad and \quad} (\rho \circ \xi \circ \eta)^{-1}(r) > r \mbox{\quad for all\ } r > 0.%
\end{equation}
Let us assume to the contrary that $s_{|\IC^2} \geq \varphi_1(r)\unit$. We want to show that this implies%
\begin{equation*}
  (\id + \rho) \circ \Gamma_i(s_{|\IC^2}) \geq \Gamma_i(s) \mbox{\quad for all\ } i \in \IC^2,%
\end{equation*}
which in turn implies $\Gamma_{\rho}(s_{|\IC^2}) \geq s_{|\IC^2}$ in contradiction to (a). We fix $i \in \IC^2$ and define%
\begin{equation*}
  a_i := [\gamma_{ij}(s_j)]_{j \in \IC_i \cap \IC^2}.%
\end{equation*}
We can assume without loss of generality that $a_i \neq 0$. Indeed, $a_i = 0$ is only possible if $\IC_i \cap \IC^2 = \emptyset$, implying%
\begin{align*}
  s_i &\leq \max\{r,\mu_i([\gamma_{ij}(r)]_{j \in \IC_i})\} \leq \max\left\{r, \gamma \circ  \left(\max_{ji \in E(\GC)}\gamma_{ij}\right)(r) \right\} \stackrel{\eqref{eq_finite_general_wlogs}}{<} \varphi_1(r)%
\end{align*}
in contradiction to our assumption hat $s_{|\IC^2} \geq \varphi_1(r)\unit$. The inequality to be shown can then be written as%
\begin{equation*}
  \rho(\mu_i(a_i)) \geq \mu_i([\gamma_{ij}(r)]_{j \in \IC_i \cap \IC^1} + a_i) - \mu_i(a_i).%
\end{equation*}
Using \eqref{eq_mafs_uniform_cont}, we see that a sufficient condition for this inequality to hold is%
\begin{equation*}
  \rho(\mu_i(a_i)) \geq \gamma( \|[\gamma_{ij}(r)]_{j \in \IC_i \cap \IC^1}\| ).%
\end{equation*}
This last inequality follows from the definition of $\varphi_1$, since%
\begin{align*}
  \rho(\mu_i(a_i)) &\geq \rho(\xi(\|a_i\|)) \geq \rho \circ \xi \circ \eta \circ \varphi_1(r) \\
	                 &= \gamma \circ \left(\max_{ji \in E(\GC)}\gamma_{ij}\right)(r) \geq \gamma( \|[\gamma_{ij}(r)]_{j \in \IC_i \cap \IC^1}\| ).%
\end{align*}
Hence, there exists at least one index $i \in \IC^2$ such that $s_i  < \varphi_1(r)$. We can split the index set $\IC$ again into%
\begin{equation*}
  \IC^1 := \left\{ i \in \IC : s_i \leq \max\{r,\varphi_1(r)\} \right\},\quad \IC^2 := \IC \backslash \IC^1,%
\end{equation*}
observing that the new $\IC^1$ contains at least one more element than the old one. With the same arguments as before, we can define another $\KC_{\infty}$-function $\varphi_2$ such that $s_i < \varphi_2(r)$ must hold for at least one $i \in \IC^2$. In this way, we can continue inductively until $\IC^2 = \emptyset$ (using that there are only finitely many indices). We thus find a $\KC_{\infty}$-function $\varphi$ such that $\|s\| \leq \varphi(r)$. Note that we have proved that $\Gamma$ satisfies the $\oplus$-MBI property (and not $\Gamma_{\rho}$). By choosing a uniformly continuous $\KC_{\infty}$-function $\rho' < \rho$, the proof also works for $\Gamma_{\rho'}$ in place of $\Gamma$.
\end{proof}

\begin{remark}
In the literature, we find a very similar result, see \cite[Thm.~C.41]{Mir}. However, here the $\oplus$-MBI property does not appear, but instead the stronger MBI property (see Remark \ref{rem_uniform_sgc}). The proof of this result uses different methods, however, and in particular uses the subadditivity of MAFs for the path construction.
\end{remark}

\begin{remark}
Due to Theorem \ref{thm_necessary_conditions}, we know that the statements in Theorem \ref{thm_finite_general_case} are also equivalent to the existence of a path satisfying properties (i), (ii) and (iv) of a $C^0$-path of strict decay, but not necessarily being monotone.
\end{remark}

\section{Some open questions}\label{sec_summary}

The main results of the paper lead to the following open questions and topics for future research:%
\begin{enumerate}
\item[1.] The ultimate goal, to obtain a complete characterization of the existence of a path of strict decay in terms of no-joint-increase conditions, was only obtained in special cases. One major step towards this goal would be to clarify when the weak$^*$-convergence of $\bsign{\Gamma}_{r\unit}^n(s)$ to $\sigma^*(r)$ for $s \geq \sigma^*(r)$ can be strengthened to norm-convergence. In particular, it must be understood if this follows from UGAS of $\Sigma(\Gamma_{\rho})$ or if stronger conditions are necessary.%
\item[2.] The relation between UGAS and the $\oplus$-MBI property should be further clarified. A related problem is the unification of the two uniform NJI conditions. From Theorem \ref{thm_finite_general_case}, we can deduce that in the classical finite case with uniformly continuous MAFs, they are equivalent. From Proposition \ref{prop_finite_case_uniform_NJI}, we further know that in the general finite case SU-NJI implies TU-NJI. A complete characterization of UGAS in terms of NJI conditions is also missing; in particular, it is unclear if UGAS implies the cofinality of the decay set.%
\item[3.] A further topic for future research concerns the question of \emph{scalable ISS} \cite{BKn}. Here, we ask under which conditions stability measures are preserved when a network is enlarged. If we can cover a large-but-finite network that may grow in time with one infinite network which is ISS, we should be able to say that stability measures are preserved when further subsystems are added. The stability measures are related to the $\KC_{\infty}$-functions appearing in the definition of a path of strict decay, so it is important to take a closer look at how these functions are related to the fundamental components of the network. A first step in this direction is Corollary \ref{cor_finite_subnetworks}.%
\item[4.] Besides answering the open questions within the setup of this paper, one can ask if this setup can be generalized. For instance, one can ask if it makes sense to build a similar theory on $\ell^p$-spaces for finite $p$, if one can relax the assumption that the sets $\IC_i$ all have to be finite, or even if the countability of the index set $\IC$ can be dropped. 
\end{enumerate}

\appendix

\section{Some technical lemmas}\label{sec_appendix}

The following three lemmas describe properties of the weak$^*$-topology of $\ell^{\infty}(\IC)$.%

\begin{lemma}\label{lem_weakstar1}
A sequence $(s^n)_{n\in\N}$ in $\ell^{\infty}(\IC)$ converges in the weak$^*$-topology if and only if it is norm-bounded and componentwise convergent.
\end{lemma}

\begin{proof}
As a general fact, the sequence $(s^n)$ converges to some $s \in \ell^{\infty}(\IC)$ if and only if $s^n(t)$ converges to $s(t)$ for every $t \in \ell^1(\IC)$. By \eqref{eq_dual_action}, this means%
\begin{equation*}
  \lim_{n \rightarrow \infty} \sum_{i \in \IC} s^n_it_i = \sum_{i \in \IC} s_it_i \mbox{\quad for all\ } t \in \ell^1(\IC).%
\end{equation*}
If we let $t$ range through the standard basis of $\ell^1(\IC)$, consisting of the unit vectors $e^i$, $i \in \IC$, this implies%
\begin{equation*}
  \lim_{n \rightarrow \infty} s^n_i = s_i\mbox{\quad for all\ } i \in \IC.%
\end{equation*}
Hence, we have shown that weak$^*$ convergence implies componentwise convergence. Now, we show that $(s^n)$ is norm-bounded. By \cite[Cor.~2.6.8]{Meg}, a set $A \subset \ell^{\infty}(\IC)$ is norm-bounded if and only if the set $\{ s(t) : s \in A \}$ is bounded in $\R$ for each $t \in \ell^1(\IC)$. Applying this to $A := \{ s^n : n \in \N \}$ shows norm-boundedness, since $s^n(t)$ is convergent (and thus bounded) for each $t$ by assumption.%

Conversely, let $(s^n)_{n\in\N}$ be a sequence in $\ell^{\infty}(\IC)$ which is norm-bounded, say $\|s^n\| \leq B$ for all $n$, and converges to some $s \in \ell^{\infty}(\IC)$ componentwise. Pick $t \in \ell^1(\IC)$ and $\ep > 0$. Obviously, we can assume that $t \neq 0$. Choose a subset $\JC \subset \IC$ such that $\IC \backslash \JC$ is finite and%
\begin{equation*}
  \sum_{i \in \JC} |t_i| < \frac{\ep}{4B}.%
\end{equation*}
The existence of such $\JC$ immediately follows from $t \in \ell^1(\IC)$. Now, we can choose $N$ large enough so that $|s^n_i - s_i| < \ep / (2\|t\|_1)$ for $i \in \IC \backslash \JC$ and $n \geq N$. The choice of $\JC$ and $N$ then implies that for all $n \geq N$%
\begin{align*}
  \Bigl|\sum_{i \in \IC}(s^n_i - s_i)t_i\Bigr| &\leq \sum_{i \in \IC}|s^n_i - s_i| \cdot |t_i| \\
	                                             &= \sum_{i \in \IC \backslash \JC} |s^n_i - s_i| \cdot |t_i| + \sum_{i \in \JC} |s^n_i - s_i| \cdot |t_i| \\
																					     &\leq \frac{\ep}{2\|t\|_1} \sum_{i \in \IC} |t_i| + \sum_{i \in \JC} (|s^n_i| + |s_i|) |t_i| \leq 
																							\frac{\ep}{2} + \frac{\ep}{2} = \ep.%
\end{align*}
This shows that $(s^n)$ converges to $s$ in the weak$^*$-topology.
\end{proof}

\begin{lemma}\label{lem_weakstar2}
Every norm-bounded sequence $(s^n)_{n\in\N}$ in $\ell^{\infty}(\IC)$ has a subsequence which is componentwise convergent.
\end{lemma}

\begin{proof}
By the theorem of Banach-Alaoglu and the assumption of norm-boundedness, $(s^n)_{n\in\N}$ is contained in a weakly$^*$ compact subset of $\ell^{\infty}(\IC)$. By \cite[Cor.~2.6.20]{Meg}, this subset is metrizable. Hence, its compactness is equivalent to its sequential compactness. The rest follows from Lemma \ref{lem_weakstar1}. 
\end{proof}

\begin{lemma}\label{lem_weakstar3}
In $\ell^{\infty}(\IC)$, every order interval of the form $[a,b]$ with $a \leq b$ is weakly$^*$ sequentially compact.
\end{lemma}

\begin{proof}
Let $(s^n)$ be a sequence in $[a,b]$. Then clearly, $(s^n)$ is norm-bounded. By Lemma \ref{lem_weakstar2}, there exists a subsequence $(s^{n_m})$ which is componentwise convergent. Let $s$ denote its limit. The inequalities $a_i \leq s^{n_m}_i \leq b_i$ which hold for all $i$ and $m$, clearly carry over to the limit when $m$ goes to infinity, so $s \in [a,b]$.
\end{proof}

The proof of the following simple lemma is left to the reader.%

\begin{lemma}\label{lem_maxdiff}
Let $a \geq c \geq 0$ and $b \geq d \geq 0$ for vectors $a,b,c,d \in \ell^{\infty}_+(\IC)$. Then%
\begin{equation*}
  a \oplus b - c \oplus d \leq (a - c) \oplus (b - d).%
\end{equation*}
\end{lemma}

\begin{lemma}{\cite[Lem.~A.6]{Mir}}\label{lem_kinfty}
For any $\rho \in \KC_{\infty}$, the following statements hold:%
\begin{enumerate}
\item[(a)] There exists $\eta \in \KC_{\infty}$ such that $(\id + \rho)^{-1} = \id - \eta$.%
\item[(b)] There are $\rho_1,\rho_2 \in \KC_{\infty}$ such that $\id + \rho = (\id + \rho_1) \circ (\id + \rho_2)$.%
\end{enumerate}
\end{lemma}

\begin{lemma}{\cite[Prop.~A.16]{Mir}}\label{lem_mir_A16}
Let $z:\R_+ \rightarrow \R_+$ be non-decreasing, continuous at zero and satisfy $z(0) = 0$, $z(r) > 0$ for all $r > 0$. Then there exist $z_1,z_2 \in \KC$, smooth on $(0,\infty)$, such that $z_1(r) \leq z(r) \leq z_2(r)$ for all $r \geq 0$. If $z(r) \rightarrow \infty$ as $r \rightarrow \infty$, then $z_1,z_2$ can be chosen in $\KC_{\infty}$.
\end{lemma}

\begin{lemma}\label{lem_ugas}
Let $T:\ell^{\infty}_+ \rightarrow \ell^{\infty}_+$ be a monotone operator with $T(0) = 0$. If the system $\Sigma(T)$ is UGS and GATT, then it is UGAS.
\end{lemma}

\begin{proof}
Fix $r,\ep > 0$ and choose $n$ such that $\|T^k(r\unit)\| \leq \ep$ for all $k \geq n$. Then $\|s\| \leq r$ implies $s \leq r\unit$, and thus $T^k(s) \leq T^k(r\unit) \leq \ep\unit$ for all $k \geq n$. Hence, $\Sigma(T)$ is UGATT. By \cite[Thm.~4.2]{Mi2}, UGS and UGATT together imply UGAS.
\end{proof}

\begin{lemma}\label{lem_decay_set}
Let $T$ be a monotone operator on $\ell^{\infty}_+(\IC)$. If $s \in \Psi(T)$ for some $s \in \ell^{\infty}_+(\IC)$, then $[T(s),s] \subset \Psi(T)$.
\end{lemma}

\begin{proof}
If $T(s) \leq s' \leq s$, then monotonicity implies $T(s') \leq T(s) \leq s'$.
\end{proof}

\end{document}